# SOME RESULTS ON THE SOLUTIONS OF CAPUTO FRACTIONAL LINEAR TIME-INVARIANT SYSTEMS OF ANY ORDER WITH INTERNAL POINT DELAYS

by M. DE LA SEN

Institute for Research and Development of Processes. Faculty of Science and Technology

University of Basque Country. Campus of Leioa. Aptdo. 544- Bilbao – SPAIN

**Abstract**. This paper is devoted to the investigation of the nonnegative solutions and the stability and asymptotic properties of the solutions of fractional differential dynamic systems involving delayed dynamics with point delays. The obtained results are independent of the sizes of the delays.

1. **Introduction**

The theory of fractional calculus is basically concerned with the calculus of integrals and derivatives of any arbitrary real or complex orders. In this sense, it may be considered as a generalization of classical calculus which is included in the theory as a particular case. The former ideas have been stated about three hundred years ago but the main mathematical developments and applications of fractional calculus have been of increasing interest from the seventies. There is a good compendium of the state of the art of the subject and the main related existing mathematical results with examples and case studies in [1]. There are a lot of results concerning the exact and approximate solutions of fractional differential equations of Riemann – Liouville and Caputo types, [1-4], fractional derivatives involving products of polynomials, [5-6], fractional derivatives and fractional powers of operators, [7-9], boundary value problems concerning fractional calculus (see, for instance, [1], [10]) etc. There is also an increasing interest in the recent mathematical literature in the characterization of dynamic fractional differential systems oriented towards several fields of Science like Physics, Chemistry or Control Theory because it is a powerful tool for later applications in all fields requiring support via ordinary, partial derivatives and functional differential equations. Perhaps the reason of interest of fractional calculus is that, the numerical value of the fraction parameter allows a closer characterization of eventual uncertainties present in the dynamic model. We can find, in particular, a lot of literature concerned with the development of Lagrangian and Hamiltonian formulations where the motion integrals are calculated though fractional calculus and also in related investigations concerned dynamic and damped and diffusive systems [11-17] as well as the characterization of impulsive responses or its use in Applied Optics related, for instance, to the formalism of fractional derivative Fourier plane filters (see, for instance, [16-18]) and Finance [19]. Fractional calculus is also of interest in Control Theory concerning for instance, heat transfer, lossless transmission lines, the use of discretizing devices supported by fractional calculus, etc. (see, for instance [20-22]). In particular, there are several recent applications of fractional calculus in the fields of filter design, circuit theory and robotics, [21-22], and signal processing, [17]. Fortunately, there is an increasing mathematical literature currently available on fractional differ-integral calculus which can formally support successfully the investigations in other related disciplines.

This paper is concerned with the investigation of the solutions of time-invariant fractional differential dynamic systems, [23-24], involving point delays what leads to a formalism of a class of functional differential equations, [25-31]. Functional equations involving point delays are a crucial mathematical tool to investigate real process where delays appear in a natural way like, for instance, transportation problems, war and peace problems or biological and medical processes. The main interest of this paper is concerned with the positivity and stability of solutions independent of the sizes of the delays and also being independent of eventual coincidence of some values of delays if those ones are, in particular, multiple related to the associate matrices of dynamics. Most of the results are centred in characterizations via Caputo fractional differentiation although some extensions are presented concerned with the classical Riemann- Liouville differ-integration. It is proved that the existence nonnegative solutions independent of the sizes of the delays and the stability properties of linear time-invariant fractional dynamic differential systems subject to point delays may be characterized with sets of precise mathematical results.

**1.1 Notation**

$Z$, $R$ and $C$ are the sets of integer, real and complex numbers, $Z_+$ and $R_+$ are the positive integer and real numbers, and

$Z_{0+} := Z_+ \cup \{0\}$; $R_{0+} := R_+ \cup \{0\}$; $C_+ := \{z \in C : Re\, z > 0\}$; $C_{0+} := \{z \in C : Re\, z \geq 0\}$

$\bar{n} := \{1, 2, ..., n\}$

The following notation is used to characterize different levels of positivity of matrices:

$R_{0+}^{n \times m} := \{M = (m_{ij}) \in R^{n \times m} : m_{ij} \geq 0; \forall (i, j) \in \bar{n} \times \bar{m}\}$ is the set of all $n \times m$ real matrices of nonnegative entries. If $M \in R^{n \times m}$ then $M \geq 0$ is used as a simpler notation for $M \in R_{0+}^{n \times m}$.

$R_+^{n \times m} := \{0 \neq M = (m_{ij}) \in R^{n \times m} : m_{ij} \geq 0; \forall (i, j) \in \bar{n} \times \bar{m}\}$ is the set of all nonzero $n \times m$ real matrices of nonnegative entries (i.e. at least one of their entries is positive). If $M \in R^{n \times m}$ then $M > 0$ is used as a simpler notation for $M \in R_+^{n \times m}$.

$R_{++}^{n \times m} := \{M = (M_{ij}) \in R^{n \times m} : M_{ij} > 0; \forall (i, j) \in \bar{n} \times \bar{m}\}$ is the set of all $n \times m$ real matrices of positive entries. If $M \in R^{n \times m}$ then $M >> 0$ is used as a simpler notation for $M \in R_{++}^{n \times m}$. The superscript T denotes the transpose, $M_i^T$ and $M_j$ are, respectively, the i-th row and the j-th column of the matrix M. A close notation to characterize the positivity of vectors is the following:

$R_{0+}^n := \{v = (v_1, v_2, ..., v_n)^T \in R^n : v_i \geq 0; \forall i \in \bar{n}\}$ is the set of all $n$ real vectors of nonnegative components. If $v \in R^n$ then $v \geq 0$ is used as a simpler notation for $v \in R_{0+}^n$.

$R_+^n := \{0 \neq v = (v_1, v_2, ..., v_n)^T \in R^n : v_i \geq 0; \forall i \in \bar{n}\}$ is the set of all $n$ real nonzero vectors of nonnegative components (i.e. at least one component is positive). If $v \in R^n$ then $v > 0$ is used as a simpler notation for $v \in R_+^n$.

$\boldsymbol{R}_{++}^{n} := \{v = (v_1, v_2, ..., v_n)^T \in \boldsymbol{R}^n : v_i > 0; \forall i \in \bar{n}\}$ is the set of all $n$ real vectors of positive components. If $v \in \boldsymbol{R}^n$ then $v >> 0$ is used as a simpler notation for $v \in \boldsymbol{R}_{++}^{n}$.

$M = (M_{ij}) \in \boldsymbol{R}^{n \times n}$ is a Metzler matrix if $M_{ij} \geq 0; \forall (i, j \neq i) \in \bar{n} \times \bar{n}$. $MR^{n \times n}$ is the set of Metzler matrices of order n.

The maximum real eigenvalue, if any, of a real matrix M, is denoted by $\lambda_{max}(M)$

## 2. Some background on fractional differential systems

Assume that $f : [a,b] \to \boldsymbol{C}^n$ for some real interval $[a,b] \subset \boldsymbol{R}$ satisfies $f \in C^{k-2}((a,b), \boldsymbol{R}^n)$ and, furthermore, $d^{k-1} f(\tau)/d\tau^{k-1}$ exists everywhere in $[a,b]$ for $k = [Re\,\alpha] + 1$ for some $\alpha \in \boldsymbol{C}_{0+}$. Then, the Riemann-Liouville left-sided fractional derivative $^{RL}D_{a+}^{\alpha} f$ of order $\alpha \in \boldsymbol{C}_{0+}$ of the vector function $f$ in $[a,b]$ is point-wise defined as:

$$\left(^{RL}D_{a+}^{\alpha} f\right)(t) := \frac{1}{\Gamma(k-\alpha)} \left(\frac{d^k}{dt^k} \int_a^t \frac{f(\tau)}{(t-\tau)^{\alpha+1-k}} d\tau\right) ; t \in [a,b] \qquad (2.1)$$

where $k = [Re\,\alpha] + 1$ and $\Gamma : \boldsymbol{C} \setminus \boldsymbol{Z}_{0-} \to \boldsymbol{C}$, where $\boldsymbol{Z}_{0-} := \{n \in \boldsymbol{Z} : n \leq 0\}$, is the $\Gamma$-function defined by $\Gamma(z) := \int_0^\infty \tau^{z-1} e^{-\tau} d\tau ; z \in \boldsymbol{C} \setminus \boldsymbol{Z}_{0-}$. If $f \in C^{k-1}((a,b), \boldsymbol{R}^n)$ and, furthermore, $f^{(k)}(\tau) \equiv d^k f(\tau)/d\tau^k$ exists everywhere in $[a,b]$, then the Caputo left-sided fractional derivative $^{C}D_{a+}^{\alpha} f$ of order $\alpha \in \boldsymbol{C}_{0+}$ of the vector function $f$ in $[a,b]$ is point-wise defined as:

$$\left(^{C}D_{a+}^{\alpha} f\right)(t) = \frac{1}{\Gamma(k-\alpha)} \int_a^t \frac{f^{(k)}(\tau)}{(t-\tau)^{\alpha+1-k}} d\tau ; t \in [a,b] \qquad (2.2)$$

where $k = [Re\,\alpha] + 1$ if $\alpha \notin \boldsymbol{Z}_{0+}$ and $k = \alpha$ if $\alpha \in \boldsymbol{Z}_{0+}$. The following relationship between both fractional derivatives holds provided that they exist (i.e. if $f : [a,b] \to \boldsymbol{C}^n$ possesses Caputo left-sided fractional derivative in $[a,b]$), [1]:

$$\left(^{C}D_{a+}^{\alpha} f\right)(t) = \,^{RL}D_{a+}^{\alpha} \left[f(\tau) - \sum_{j=0}^{k-1} \frac{f^{(j)}(a)(\tau-a)^j}{j!}\right](t); t \in [a,b] \qquad (2.3)$$

Since $Re\,\alpha \leq k$, the above formula relating both fractional derivatives proves the existence of the Caputo left-sided fractional derivative in $[a,b]$ if that of Riemann – Liouville exists in $[a,b]$.

## 3. Solution of a fractional differential dynamic system of any order $\alpha$ with internal point delays

Consider the linear and time-invariant differential functional Caputo fractional differential system of order $\alpha$:

$$\left(^{C}D_{0+}^{\alpha} x\right)(t) = \sum_{i=0}^{p} A_i x(t - h_i) + Bu(t) \qquad (3.1)$$

with $k-1 < \alpha (\in \mathbf{R}_+) \leq k$; $k-1, k \in \mathbf{Z}_{0+}$, $0 = h_0 < h_1 < h_2 < ... < h_p = h < \infty$ being distinct constant delays, $A_0, A_i \in \mathbf{R}^{n \times n}$ ($i \in \bar{p} := \{1, 2, ...., p\}$), are the matrices of dynamics for each delay $h_i, i \in \bar{p} \cup \{0\}$, $B \in \mathbf{R}^{n \times m}$ is the control matrix. The initial condition is given by k n-real vector functions $\varphi_j : [-h, 0] \to \mathbf{R}^n$, with $j \in \overline{k-1} \cup \{0\}$, which are absolutely continuous except eventually in a set of zero measure of $[-h, 0] \subset \mathbf{R}$ of bounded discontinuities with $\varphi_j(0) = x_j(0) = x^{(j)}(0) = x_{j0}$. The function vector $u : \mathbf{R}_{0+} \to \mathbf{R}^m$ is any given bounded piecewise continuous control function. The following result is concerned with the unique solution on $\mathbf{R}_{0+}$ of the above differential fractional system (3.1). The proof follows directly from a parallel existing result from the background literature on fractional differential systems by grouping all the additive forcing terms of (3.1) in a unique one (see, for instance [1], Eqs. (1.8.17), (3.1.34)-(3.1.49), with $f(t) \equiv \sum_{i=1}^{p} A_i x(t - h_i) + Bu(t)$).

**Theorem 3.1**. The linear and time-invariant differential functional fractional differential system (3.1) of any order $\alpha \in C_{0+}$ has a unique solution on $\mathbf{R}_{0+}$ for each given set of initial functions $\varphi_j : [-h, 0] \to \mathbf{R}^n$, $j \in \overline{k-1} \cup \{0\}$ being absolutely continuous except eventually in a set of zero measure of $[-h, 0] \subset \mathbf{R}$ of bounded discontinuities with $\varphi_j(0) = x_j(0) = x^{(j)}(0) = x_{j0}$; $j \in \overline{k-1} \cup \{0\}$ and each given control $u : \mathbf{R}_{0+} \to \mathbf{R}^m$ being a bounded piecewise continuous control function. Such a solution is given by:

$$x_\alpha(t) = \sum_{j=0}^{k-1} \left( \Phi_{\alpha j 0}(t) x_{j0} + \sum_{i=1}^{p} \int_0^{h_i} \Phi_\alpha(t - \tau) A_i \varphi_j(\tau - h_i) d\tau \right)$$
$$+ \sum_{i=1}^{p} \int_{h_i}^{t} \Phi_\alpha(t - \tau) A_i x_\alpha(\tau - h_i) d\tau + \int_0^t \Phi_\alpha(t - \tau) Bu(\tau) d\tau \; ; \; t \in \mathbf{R}_{0+} \qquad (3.2)$$

with $k = [Re \, \alpha] + 1$ if $\alpha \notin \mathbf{Z}_+$ and $k = \alpha$ if $\alpha \in \mathbf{Z}_+$, and

$$\Phi_{\alpha j 0}(t) := t^j E_{\alpha, j+1}(A_0 t^\alpha) \; ; \; \Phi_\alpha(t) := t^{\alpha-1} E_{\alpha, \alpha}(A_0 t^\alpha) \qquad (3.3)$$

$$E_{\alpha, j}(A_0 t^\alpha) := \sum_{\ell=0}^{\infty} \frac{(A_0 t^\alpha)^\ell}{\Gamma(\alpha \ell + j)} \; ; \; j \in \overline{k-1} \cup \{0, \alpha\} \qquad (3.4)$$

for $t \geq 0$ and $\Phi_{\alpha 0}(t) = \Phi_\alpha(t) = 0$ for $t < 0$, where $E_{\alpha, j}(A_0 t^\alpha)$ are the Mittag-Leffler functions. □□

Now consider that the right-hand-side of (3.1) is the evaluation of a Riemann-Liouville fractional differential system of the same order $\alpha$ as follows:

$$\left( {}^{RL}D_{0+}^\alpha x \right)(t) = \sum_{i=0}^{p} A_i x(t - h_i) + Bu(t) \qquad (3.5)$$

under the same functions of initial conditions as those of (3.1). Through the formula (2.3) relating Caputo and Riemann-Liouville left-sided fractional derivatives of the same order $\alpha$, one gets:

$$\left({}^C D_{0+}^{\alpha} x\right)(t) = \sum_{i=0}^{p} A_i x(t-h_i) + Bu(t) - {}^{RL}D_{0+}^{\alpha}\left(\sum_{j=0}^{k-1} \frac{x_{j0} \tau^j}{j!}\right)(t) \tag{3.6}$$

Since the Caputo left-sided fractional derivative and the Riemann-Liouville fractional integral of order $\alpha \in \hat{C}_+ := \{Z_+ \cup \{z \in C_+ : Re\ z \notin Z_+\}\}$ are inverse operators (what is not the case if $\alpha \notin \hat{C}_+$),(see [1], Lemma 2.21(a)), one gets from (3.6), (2.3) and (3.2) if $\alpha \in \hat{C}$ the subsequent result for the fractional differential system (3.5) on $R_{0+}$:

**Corollary 3.2**. If (3.5) of any order $\alpha \in \hat{C}_+$ is replaced with (3.1) under the same initial conditions then its unique solution on $R_{0+}$ is given by:

$$x_\alpha(t) = \sum_{j=0}^{k-1}\left(\left(\Phi_{\alpha j0}(t) - \frac{t^j}{j!}I_n\right)x_{j0} + \sum_{i=1}^{p}\int_0^{h_i}\Phi_\alpha(t-\tau)A_i\varphi_j(\tau-h_i)d\tau\right)$$

$$+ \sum_{i=1}^{p}\int_{h_i}^{t}\Phi_\alpha(t-\tau)A_i x_\alpha(\tau-h_i)d\tau + \int_0^t \Phi_\alpha(t-\tau)Bu(\tau)d\tau\ ;\ \alpha \notin Z_{0+}\ ;\ t \in R_{0+}$$

$$= \sum_{j=0}^{k-1}\left(\left(E_{\alpha,j+1}(t) - \frac{1}{j!}I_n\right)t^j x_{j0} + \sum_{i=1}^{p}\int_0^{h_i}\Phi_\alpha(t-\tau)A_i\varphi_j(\tau-h_i)d\tau\right)$$

$$+ \sum_{i=1}^{p}\int_{h_i}^{t}\Phi_\alpha(t-\tau)A_i x_\alpha(\tau-h_i)d\tau + \int_0^t \Phi_\alpha(t-\tau)Bu(\tau)d\tau\ ;\ \alpha \notin Z_{0+}\ ;\ t \in R_{0+} \tag{3.7}$$

with $k = [Re\ \alpha]+1$ if $\alpha \notin Z_+$ and $k = \alpha$ if $\alpha \in Z_+$. □

Another mild evolution operator can be considered to construct the unique solution of (3.1) by considering the control effort as the unique forcing term of (3.1) and the functions of initial conditions as forcing terms. See the corresponding expressions obtainable from [1], Eqs. (1.8.17), (3.1.34)-(3.1.49), with the identity $f(t) \equiv Bu(t)$) and the evolution operator defined in [2-3] for the standard (non-fractional differential system), i.e. $\alpha = 1$ in (3.1). Thus, another equivalent expression for the unique solution of the Caputo fractional differential system of order $\alpha$ is given in the subsequent result:

**Theorem 3.3**. The solution of (3.1) given in Theorem 3.1 is equivalently rewritten as follows:

$$x_\alpha(t) = \sum_{j=0}^{k-1}\left(\Psi_{\alpha j0}(t)x_{j0} + \sum_{i=1}^{p}\int_0^{h_i}\Psi_{\alpha j0}(t-\tau)\varphi_j(\tau-h_i)d\tau\right) + \int_0^t \Psi_\alpha(t-\tau)Bu(\tau)d\tau \tag{3.8}$$

for $t \in R_{0+}$, any $\alpha \in C_+$ with $k = [Re\ \alpha]+1$ if $\alpha \notin Z_+$ and $k = \alpha$ if $\alpha \in Z_+$; and

$$\Psi_{\alpha j0}(t) := t^j E_{\alpha,j+1}\left(A_0 t^\alpha\right) + \sum_{i=1}^{p}\int_0^t \tau^{\alpha-1} E_{\alpha,\alpha}\left(A_0 \tau^\alpha\right) A_i \Psi_{\alpha j0}(t-\tau-h_i)d\tau \tag{3.9}$$

$$\Psi_\alpha(t) := t^{\alpha-1} E_{\alpha,\alpha}(A_0 t^\alpha) + \sum_{i=1}^{p} \int_0^t \tau^{\alpha-1} E_{\alpha,\alpha}(A_0 \tau^\alpha) A_i \Psi_\alpha(t-\tau-h_i) d\tau \qquad (3.10)$$

for $t \geq 0$ and $\Psi_{\alpha j 0}(t) = \Psi_\alpha(t) = 0$, $j \in \overline{k-1} \cup \{0\}$ for $t \in [-h, 0)$. □

Also, the solution to the Riemann-Liouville fractional differential system (3.5) under the same initial conditions as those of (4) is given in the next result for $k = [Re\ \alpha] + 1$ if $\alpha \notin \mathbf{Z}_+$ based on (3.6):

**Corollary 3.4.** If (3.5) being of order $\alpha \in \hat{\mathbf{C}}_+$ is replaced with (3.1) under the same initial conditions then its unique solution on $\mathbf{R}_{0+}$ is given by

$$x_\alpha(t) = \sum_{j=0}^{k-1} \left[ \left( \Psi_{\alpha j 0}(t) - \frac{t^j}{j!} I_n \right) x_{j0} + \sum_{i=1}^{p} \int_0^{h_i} \Psi_{\alpha j 0}(t-\tau) \varphi_j(\tau - h_i) d\tau \right] + \int_0^t \Psi_\alpha(t-\tau) B u(\tau) d\tau$$

$$(3.11)$$

with $k = [Re\ \alpha] + 1$ if $\alpha \notin \mathbf{Z}_+$ and $k = \alpha$ if $\alpha \in \mathbf{Z}_+$ which is identical to that given in Corollary 3.2. □

Particular cases of interest of the solution of (3.1) given in Theorem 3.3 are:

1) $\alpha = k$ which yields the solution:

$$x_k(t) = \sum_{j=0}^{k-1} \left( \Psi_{kj0}(t) x_{j0} + \sum_{i=1}^{p} \int_0^{h_i} \Psi_{kj0}(t-\tau) \varphi_j(\tau - h_i) d\tau \right) + \int_0^t \Psi_k(t-\tau) B u(\tau) d\tau \qquad (3.12)$$

2) A further particular case $\alpha = k = 1$ yields the solution:

$$x_1(t) = \Psi_1(t) x_{j0} + \sum_{i=1}^{p} \int_0^{h_i} \Psi_1(t-\tau) \varphi_0(\tau - h_i) d\tau + \int_0^t \Psi_1(t-\tau) B u(\tau) d\tau \qquad (3.13)$$

since $\Psi_{100}(t) = \Psi_1(t)$, $t \in \mathbf{R}_{0+}$ which is the unique solution of $(Dx)(t) = \sum_{i=0}^{p} A_i(t-h_i) + Bu(t)$ under any almost everywhere absolutely continuous function (except eventually in some subset of zero measure of $[-h, 0]$ of bounded discontinuities) of initial conditions $\varphi \equiv \varphi_0 : [-h, 0] \to \mathbf{R}^n$. Use for this case, the less involved notations $(\Psi x)(t) = (\Psi_{100} x)(t) = (\Psi_1 x)(t)$ for the smooth evolution operator from $\mathbf{R}_{0+} \times \mathbf{R}^n$ to $\mathbf{R}^n$, and $\Phi(t) = \Phi_{100}(t) = \Phi_1(t) = e^{A_0 t}$, $t \in \mathbf{R}_+$ for the exponential matrix function $e^{A_0 t}$ from $\mathbf{R}_{0+}$ to $\mathbf{R}^{n \times n}$, which defines a $C_0$-semigroup $(e^{A_0 t}, t \in \mathbf{R}_{0+})$ of infinitesimal generator $A_0$ from $\mathbf{R}_{0+}$ to $L(\mathbf{R}^n)$. Then, the unique solution $x(t) \equiv x_1(t)$, $t \in \mathbf{R}_+$ for the given function of initial conditions is:

$$x(t) = \Phi(t) x_0 + \sum_{i=1}^{p} \int_0^{h_i} \Phi(t-\tau) A_i \varphi(\tau - h_i)$$

$$+ \sum_{i=1}^{p} \int_{h_i}^{t} \Phi(t-\tau) A_i x(\tau - h_i) d\tau + \int_0^t \Phi(t-\tau) B u(\tau) d\tau \qquad (3.14)$$

$$= \Psi(t)x_0 + \sum_{i=1}^{p} \int_{0}^{h_i} \Psi(t-\tau)\varphi(\tau - h_i)d\tau + \int_{0}^{t} \Psi(t-\tau)Bu(\tau)d\tau \ ; \ t \in \mathbf{R}_{0+} \quad (3.15)$$

and $x(t) = \varphi(t)$ for $t \in [-h, 0]$, where $\Phi(t) = e^{A_0 t}$ satisfies $\dot{\Phi}(t) = A_0 \Phi(t)$ $t \in \mathbf{R}$, and $\dot{\Psi}(t) = \sum_{i=0}^{p} A_i \Psi(t - h_i)$ with $\Psi(0) = \Phi(0) = I_n$ (the n- identity matrix) and $\Psi(t) = 0$, $t \in [-h, 0)$ which has a unique solution $\Psi(t) = e^{A_0 t} \left( I_n + \sum_{i=1}^{p} \int_{h_i}^{t} e^{-A_0 \tau} A_i \Psi(\tau - h_i) d\tau \right)$ for $t \in \mathbf{R}_{0+}$, [2-3].

A problem of interest when considering a set of p delays in $[0, h]$ is the case of potentially repeated delays, then subject to $0 = h_0 \le h_1 \le h_2 \le ... \le h_p \le h < \infty$, with q of them $\{h_{jp}, j \in \bar{q} \cup \{0\}\}$ being distinct, each being repeated $1 \le v_j \le p$ ($j \in \bar{q} \cup \{0\}$) times so that

$$0 = h_0 = h_{0p} < h_{1p} < ... < h_{qp} \le h < \infty \ ; \ \sum_{j=0}^{p} v_{jp} = p+1 \quad (3.16)$$

$$h_{jp} = h_{k+i}, \ k = \sum_{\ell=0}^{j-1} v_\ell \ ; \ \forall i \in \bar{v_j}, \ j \in \bar{q} \cup \{0\} \quad (3.17)$$

Thus, the following result holds from Theorem 3.3 by grouping the terms of the delayed dynamics corresponding to the same potentially repeated delays.

**Theorem 3.5**. The Caputo solutions to the subsequent Caputo and Riemann- Liouville fractional differential systems of order $\alpha$ with $p \ge 0$ (potentially repeated) delays and $0 \le q \le p$ distinct delays:

$$\left( {}^{C}D_{0+}^{\alpha} x \right)(t) = \sum_{i=0}^{p} A_i (t - h_i) + Bu(t) \text{ and } \left( {}^{RL}D_{0+}^{\alpha} x \right)(t) = \sum_{i=0}^{p} A_i (t - h_i) + Bu(t)$$

on $\mathbf{R}_{0+}$ for the given set of initial conditions on $[-h, 0)$ are given by:

$$x_\alpha(t) = \sum_{j=0}^{k-1} \left( \Psi_{\alpha j 0}(t) x_{j0} + \sum_{i=1}^{q} \int_{0}^{h_{ip}} \Psi_{\alpha j 0}(t-\tau) \varphi_j (\tau - h_{ip}) d\tau \right)$$

$$+ \int_{0}^{t} \Psi_\alpha(t-\tau) Bu(\tau) d\tau \quad (3.18)$$

for any $\alpha \in \mathbf{C}_+$ with $k = [Re \ \alpha] + 1$ if $\alpha \notin \mathbf{Z}_{0+}$ and $k = \alpha$ if $\alpha \in \mathbf{Z}_{0+}$; and, respectively by

$$x_\alpha(t) = \sum_{j=0}^{k-1} \left( \left( \Psi_{\alpha j 0}(t) - \frac{t^j}{j!} \right) x_{j0} + \sum_{i=1}^{q} \int_{0}^{h_{ip}} \Psi_{\alpha j 0}(t-\tau) \varphi_j (\tau - h_{ip}) d\tau \right)$$

$$+ \int_{0}^{t} \Psi_\alpha(t-\tau) Bu(\tau) d\tau \quad (3.19)$$

for any $\alpha \in \hat{\mathbf{C}}_+$ with $k = [Re \ \alpha] + 1$ if $\alpha \notin \mathbf{Z}_+$ and $k = \alpha$ if $\alpha \in \mathbf{Z}_+$, where:

$$\Psi_{\alpha j 0}(t) := E_{\alpha, j+1} \left( \left( \sum_{i=0}^{v_0 - 1} A_i \right) t^\alpha \right) + \sum_{i=1}^{q} \int_{0}^{t} \tau^{\alpha - 1} E_{\alpha, \alpha} \left( \left( \sum_{i=0}^{v_0 - 1} A_i \right) \tau^\alpha \right) \left( \sum_{\ell=1}^{v_i} A_{\left( \sum_{j=0}^{i-1} v_j + \ell \right)} \right) \Psi_{\alpha j 0}(t - \tau - h_i) d\tau$$

$$(3.20)$$

$t \geq 0$ and $\Psi_{\alpha j 0}(t) = \Psi_\alpha(t) = 0$, $j \in \overline{k-1} \cup \{0\}$ for $t \in [-h, 0)$. □

## 4. Nonnegativity of the solutions

The positivity of the solutions of (3.1) independent of the values of the delays is now investigated under initial conditions $\varphi_j : [-h, 0] \to R_{0+}^n$, $j \in \overline{k-1} \cup \{0\}$.

**Theorem 4.1**. The Caputo fractional differential system (3.1) under the delay constraint $0 = h_0 < h_1 < h_2 < ... < h_p = h < \infty$ for any given absolutely continuous functions of initial conditions $\varphi_j : [-h, 0] \to R_{0+}^n$, $j \in \overline{k-1} \cup \{0\}$ and any piecewise continuous vector function $u : R_{0+} \to R_{0+}^n$ if $k = [\alpha] + 1$ if $\alpha \notin Z_+$ and $k = \alpha \in Z_+$, $\forall t \in R_{0+}$; $\forall \alpha \in R_+$ has following properties:

(i) $\Phi_{\alpha j 0}(t)$ is nonsingular; $\forall j \in \overline{k-1} \cup \{0\}$ and $\Phi_\alpha(t) \geq 0$; $\forall t \in R_{0+}$ (if $B \in R_+^{n \times m}$ then $\Phi_\alpha(t) > 0$; $\forall t \in R_{0+}$).

(ii)

(ii.1) $A_0 \in MR^{n \times n} \Leftrightarrow \Phi(t) \equiv \Phi_{100}(t) > 0$; $\forall t \in R_{0+}$

(ii.2) $A_0 \in MR^{n \times n} \Rightarrow \Phi_{\alpha j 0}(t) > 0$; $\forall j \in \overline{k-1} \cup \{0\}$; $\forall t \in [0, \bar{t}_j)$ for some sufficiently small $\bar{t}_j \in R_+$ with $\Phi_{\alpha 00}(t) > 0$, $\forall t \in R_{0+}$ (i.e. $\bar{t}_0 = \infty$). This property holds $\forall t \in R_{0+}$ (i.e. $\bar{t}_j = \infty$; $\forall j \in \overline{k-1} \cup \{0\}$) if, in addition, either $A_0 \geq 0$ or if $A_0$ is nilpotent or if $0 < \alpha \leq k = 1$. Furthermore, there are at least n entries (one per row) of $\Phi_{\alpha j 0}(t)$ being positive; $\forall t \in R_{0+}$. Furthermore,

(iii) Any solution (3.2) to any Caputo fractional differential system (3.1) is nonnegative independent of the delays; i.e. $x_\alpha(t) \in R_{0+}^n$; $\forall t \in [-h, \bar{t}) \cap R_{0+}$ for some $\bar{t} \in R_{0+}$, for any set of delays satisfying $0 = h_0 < h_1 < h_2 < ... < h_p \leq h < \infty$ and any absolutely continuous functions of initial conditions $\varphi_j : [-h, 0] \to R_{0+}^n$, $\forall j \in \overline{k-1} \cup \{0\}$ and any piecewise continuous control $u : R_{0+} \to R_{0+}^m$, if and only if $A_0 \in MR^{n \times n}$ for $\bar{t} \in R_{0+}$ being sufficiently small. Furthermore, $x_\alpha(t) \in R_{0+}^n$; $\forall t \in [-h, 0) \cup R_{0+}$ if, in addition, either $A_0 \geq 0$ or if $A_0$ is nilpotent or if $0 < \alpha \leq k = 1$, $A_i \in R_{0+}^{n \times n}$ ($\forall i \in \bar{p}$) and $B \in R_{0+}^{n \times m}$.

**Proof**: It is now proven that $\Phi_{\alpha j 0}(t) \geq 0$; $\forall t \in R_{0+} \Rightarrow A_0 \in MR^{n \times n}$; $\forall \alpha \in R_+$ for any $j \in \overline{k-1} \cup \{0\}$.

First, note the following. If $\alpha = k = 1$ then $\Phi_{\alpha 00}(t) = E_{1,1}(A_0 t) = \Phi(t) = \sum_{\ell=0}^{\infty} \frac{A_0^\ell t^\ell}{\ell!} = e^{A_0 t} \geq 0$ if $A_0 \in MR^{n \times n}$ from the above part of the proof and also $A_0 \in MR^{n \times n} \Rightarrow \Phi(t) \geq 0$; $\forall t \in R_{0+}$. This follows by contradiction. Assume that $\Phi_{im}(t) < 0$ for some $t \in R_+$. Consider the positive differential system $\dot{x}(t) = A_0 x(t)$, $x(0) = e_j$, $A_0 \in MR^{n \times n}$ so that $x_i(t) = -|\Phi_{im}(t)| < 0$ which contradicts the system being positive. Thus, $A_0 \in MR^{n \times n} \Leftrightarrow \Phi(t) \geq 0$; $\forall t \in R_{0+}$. Furthermore, since $\Phi(t)$ is a fundamental matrix of solutions of the differential system, it is non-singular for all finite time and the above result is weakened as follows:

$A_0 \in MR^{n \times n} \Leftrightarrow (\Phi(t) = \Phi_{\alpha 00}(t) = e^{A_0 t} > 0 \wedge \Phi(t)$ is non-singular; $\forall t \in R_{0+})$. Since $\Phi(t)$ is nonsigular; $\forall t \in R_{0+}$ at least n of its entries ( one per-row) is positive. Property (i) has been proven. Now, one gets from (3.3)-(3.4):

$$\Phi_{\alpha j 0}(t) = E_{\alpha, j+1}(A_0 t^{\alpha}) = \sum_{\ell=0}^{\infty} \frac{A_0^{\ell} t^{\alpha \ell}}{\Gamma(\alpha \ell + j + 1)} = \sum_{\ell=0}^{\infty} \frac{A_0^{\ell} t^{\ell}}{\ell!} \frac{t^{(\alpha-1)\ell} \ell!}{(\alpha \ell + j)\Gamma(\alpha \ell + j)}; j \in \overline{k-1} \cup \{0\} \quad (4.1)$$

Let $e_i$ the i-th unit Euclidean vector of $R^n$ whose i-th component is 1. Then, one obtains for all $j \in \overline{k-1} \cup \{0\}$, irrespective of the value of $\alpha \in R_+$ and $k \in Z_+$ being $k = [\alpha] + 1$ if $\alpha \notin Z_+$ and $k = \alpha \in Z_+$, provided that $A_0 \in MR^{n \times n}$:

$$\left(\Phi_{\alpha j 0}(t)\right)_{im} = e_i^T \Phi_{\alpha j 0}(t) e_m = e_i^T \left( \sum_{\ell=0}^{\infty} \frac{A_0^{\ell} t^{\ell}}{\ell!} \frac{t^{(\alpha-1)\ell} \ell!}{\Gamma(\alpha \ell + j + 1)} \right) e_m$$

$$= e_i^T \left( \sum_{\ell=0}^{\infty} \frac{A_0^{\ell} t^{\ell}}{\ell!} \frac{t^{(\alpha-1)\ell} \ell!}{(\alpha \ell + j)\Gamma(\alpha \ell + j)} \right) e_m$$

$$\geq e_i^T \left( e^{A_0 t} \right) e_m \min_{0 \leq \ell \leq \infty} \left( \frac{t^{(\alpha-1)\ell} \ell!}{\Gamma(\alpha \ell + k)} \right) = e_i^T \left( \sum_{\ell=0}^{\infty} \frac{A_0^{\ell} t^{\ell}}{\ell!} \frac{t^{(\alpha-1)\ell} \ell!}{(\alpha \ell + j)\Gamma(\alpha \ell + j)} \right) e_m \quad (4.2)$$

$$= e_i^T \left( \sum_{\ell=0}^{N} \frac{A_0^{\ell} t^{\ell}}{\ell!} \right) e_m \min_{0 \leq \ell \leq N} \left( \frac{t^{(\alpha-1)\ell} \ell!}{(\alpha \ell + j)\Gamma(\alpha \ell + k)} \right)$$

$$= e_i^T \left( e^{A_0 t} \right) e_m \min_{0 \leq \ell \leq N} \left( \frac{t^{(\alpha-1)\ell} \Gamma(\ell+1)}{(\alpha \ell + j)\Gamma(\alpha \ell + k)} \right) \geq e_i^T \left( e^{A_0 t} \right) e_m \min_{0 \leq \ell \leq N} \left( \frac{t^{(\alpha-1)\ell} \Gamma(\ell+1)}{(\alpha \ell + j)\Gamma(\alpha \ell + k)} \right) \quad (4.3)$$

$\forall (i,m) \in \overline{n} \times \overline{n}$, $\forall t \in R_{0+}$ since $A_0 \in MR^{n \times n} \Leftrightarrow \Phi(t) \geq 0$; $\forall t \in R_{0+}$, for some $N(\leq \infty) \in Z_{0+}$ and N is finite if and only if $A_0$ is nilpotent (of degree N) . Eq. (4.3) implies that $e_i^T \left( e^{A_0 t} \right) e_m \geq 0$ and then $\left(\Phi_{\alpha j 0}(t)\right)_{im} \geq 0$, $\forall (i,m) \in \overline{n} \times \overline{n}$ in the following cases:

(a) $N \leq \infty$, $t \in [0, \bar{t})$, since $\cap R_{0+}$ and some sufficiently small $\bar{t} \in R_+$, since $A_0 \in MR^{n \times n} \Rightarrow I + A_0 t > 0$; $\forall t \in R_{0+}$ and

$\Phi_{\alpha j 0}(t) = \sum_{\ell=0}^{N} \frac{A_0^{\ell} t^{\ell}}{\ell!} \frac{t^{(\alpha-1)\ell} \ell!}{(\alpha \ell + j)\Gamma(\alpha \ell + j)} = I + A_0 t + o(t) > 0$ for some sufficiently small $t \in R_{0+}$ for any $\alpha \in R_+$.

(b) $N \leq \infty$ and $A_0 \geq 0$ since $\dfrac{t^{(\alpha-1)\ell}\ell!}{(\alpha\ell+j)\Gamma(\alpha\ell+j)} \geq 0$ ; $\forall \ell \in \mathbf{Z}_{0+}$, $\forall j \in \overline{k-1} \cup \{0\}$ for any $\alpha \in \mathbf{R}_+$. It follows from inspection of (4.2) since $e_i^T(e^{A_0 t})e_m \geq 0$ $\forall (i,m) \in \overline{n} \times \overline{n}$, since $A_0 \in M\mathbf{R}^{n \times n}$. This implies $(\Phi_{\alpha j 0}(t))_{im} \geq 0$ ; $\forall t \in \mathbf{R}_{0+}$.

(c) $N < \infty \Rightarrow \min\limits_{0 \leq \ell \leq N}\left(\dfrac{t^{(\alpha-1)\ell}}{\Gamma(\alpha\ell+k)}\right) > 0$ ; $\forall j \in \overline{k-1} \cup \{0\}$ for any $\alpha \in \mathbf{R}_+$ so that $e_i^T(e^{A_0 t})e_m \geq 0$, $\forall (i,m) \in \overline{n} \times \overline{n}$, since $A_0 \in M\mathbf{R}^{n \times n}$, irrespectively of $A_0 \geq 0$ or not, what follows from (4.3). This implies $(\Phi_{\alpha j 0}(t))_{im} \geq 0$ ; $\forall t \in \mathbf{R}_{0+}$.

(d) $N \leq \infty$, $0 < \alpha \leq k = 1$. Then, $j=0$ so that

$$\dfrac{t^{(\alpha-1)\ell}\ell!}{\Gamma(\alpha\ell+j+1)} = \dfrac{t^{(\alpha-1)\ell}\ell!}{\Gamma(\alpha\ell+1)} = \dfrac{t^{(\alpha-1)\ell}\ell!}{\alpha\ell\Gamma(\alpha\ell)} = \dfrac{t^{(\alpha-1)\ell}(\ell-1)!}{\alpha\Gamma(\alpha\ell)} = \dfrac{\Gamma(\ell)}{\alpha t^{(1-\alpha)\ell}\Gamma(\alpha\ell)} \geq \dfrac{1}{\alpha t^{(1-\alpha)\ell}}$$

; $\forall \ell \in \mathbf{Z}_{0+}$ ; $\forall t \in \mathbf{R}_{0+}$ since $0 < \alpha \leq 1$ implies

$$\Gamma(\alpha\ell) = \int_0^\infty \tau^{\alpha\ell-1} e^{-\tau} d\tau \leq \Gamma(\ell) = \int_0^\infty \tau^{\ell-1} e^{-\tau} d\tau, \forall \ell \in \mathbf{Z}_{0+}$$

As a result, $\Phi_{\alpha 00}(t)$ from (4.2); $\forall t \in \mathbf{R}_{0+}$. Also, direct calculations with (3.3)-(3.4) lead to

$$\Phi_\alpha(t) := t^{\alpha-1} E_{\alpha,\alpha}(A_0 t^\alpha) = \sum_{\ell=0}^\infty \dfrac{t^{\alpha-1} A_0^\ell t^{\alpha\ell}}{\Gamma((\ell+1)\alpha)} = \sum_{\ell=0}^\infty \dfrac{A_0^\ell t^\ell}{\ell!} \dfrac{t^{(\alpha-1)(\ell+1)}\ell!}{\Gamma((\ell+1)\alpha)} \tag{4.4}$$

and similar developments to the above ones yield $(\Phi_\alpha(t))_{im} \geq 0$; $\forall (i,m) \in \overline{n} \times \overline{n}$, $\forall t \in \mathbf{R}_{0+}$ under the same conditions as above in the cases (a) to (d) for $\Phi_{\alpha 00}(t)$. On the other hand, one gets from (3.2)-(3.4) for the unforced system with point initial conditions at t=0:

$$x_\alpha(t) = \sum_{j=0}^{k-1} \Phi_{\alpha j 0}(t) x_{j0} = [\Phi_{\alpha 00}(t), \ldots, \Phi_{\alpha,k-1,0}(t)][x_{00}^T, \ldots, x_{k-1,0}^T]^T$$

which leads to $x_\alpha(t) = \Phi_{\alpha i 0}(t) x_{i0}$ by taking point initial conditions $x_{i0} \neq 0$, $x_{j0} = 0$, $(i \neq j), j \in \overline{k-1} \cup \{0\}$ so that $\Phi_{\alpha i 0}(t)$ is nonsingular for all $t \in \mathbf{R}_{0+}$ since otherwise the solution is not unique for each given set of initial conditions since any trajectory solution subject to some set of initial conditions $x_{i0} \neq 0$, $x_{j0} = 0$, would have infinitely many initial conditions, subject to identical constraint, so that such a trajectory is not unique which is a contradiction. Since this reasoning may be made for any $j \in \overline{k-1} \cup \{0\}$, $\Phi_{\alpha j 0}(t)$ is nonsingular for all $j \in \overline{k-1} \cup \{0\}$, all and, in addition, $\Phi_{\alpha j 0}(t) > 0$; $\forall j \in \overline{k-1} \cup \{0\}, \forall t \in \mathbf{R}_{0+}$ if either $A_0 \geq 0$ or if $A_0$ is nilpotent or if $0 < \alpha \leq k = 1$ or without these restricting condition within some first interval $[0, \bar{t})$. The following properties have been proven:

a) $A_0 \in M\mathbf{R}^{n \times n} \Leftrightarrow \Phi_{\alpha 00}(t) > 0$ ; $\forall t \in \mathbf{R}_{0+}$

b) $A_0 \in MR^{n \times n} \Rightarrow (\Phi_{\alpha j0}(t) > 0 \wedge \det \Phi_{\alpha j0}(t) \neq 0 \wedge \Phi_\alpha(t) \geq 0$ ; $\forall j \in \overline{k-1} \cup \{0\}$, $k = [\alpha]+1$ if $\alpha \notin Z_+$ and $k = \alpha \in Z_+$, $\forall t \in R_{0+})$

It remains to prove $\Phi_{\alpha j0}(t) > 0; \forall t \in R_{0+} \Rightarrow A_0 \in MR^{n \times n}$; $\forall j \in \overline{k-1}$, some $t \in R_{0+}$. This is equivalent to its contra-positive logic proposition. Proceed by contradiction by assuming $\exists j \in \overline{k-1}$ such that $A_0 \notin MR^{n \times n} \Rightarrow \Phi_{\alpha j0}(t) < 0$, some $t \in R_{0+}$. Note that

$$A_0 \notin MR^{n \times n} \Rightarrow e_i^T e^{A_0 t}(t) e_m = e_i^T \Phi_{\alpha 00}(t) e_m < 0, \text{ some } t \in R_{0+}, \text{ some } (i, m \neq i) \in \overline{n}$$

Then, one gets:

$$\left(\Phi_{\alpha j0}(t)\right)_{im} = e_i^T \Phi_{\alpha j0}(t) e_m \leq e_i^T \left(e^{A_0 t}\right) e_m \max_{0 \leq \ell \leq \infty}\left(\frac{t^{(\alpha-1)\ell} \ell!}{\Gamma(\alpha\ell+k)}\right) < 0$$

what contradicts $\left(-\Phi_{\alpha j0}(t)\right) > 0; \forall t \in R_{0+} \Leftarrow A_0 \notin MR^{n \times n}$; $\forall j \in \overline{k-1}$, some $t \in R_{0+}$. Thus, the proof of Properties (i) –(ii) becomes complete since the above proven property a) extends to any $j \in \overline{k-1} \cup \{0\}$ as follows:

c) $A_0 \in MR^{n \times n} \Leftrightarrow \Phi_{\alpha j0}(t) > 0$ ; $\forall j \in \overline{k-1} \cup \{0\}$, $k = [\alpha]+1$ if $\alpha \notin Z_+$ and $k = \alpha \in Z_+$, $\forall t \in R_{0+}$; $\forall \alpha \in R_+$

so that the unforced solution for any set of nonnegative point initial conditions is nonnegative for all time and, furthermore, $\varphi_i(t) \geq 0 (\forall i \in \overline{k-1} \cup \{0\}); \forall t \in [-h, 0], u(t) \in R_{0+}^n; \forall t \in R_{0+}, A_i \geq 0 (\forall i \in \overline{p})$ and $B \geq 0$; $\forall t \in R_{0+}$ implies that (3.2) is everywhere nonnegative within its definition domain. The converse is also true as it follows by contradiction arguments. If there is one entry of B or $A_i$ (some $i \in \overline{p}$) which is negative, or if $A_0 \notin MR^{n \times n}$, it can always be found a control $u(t) \in R_{0+}^n$ of sufficiently large norm along a given time interval such that some component of the solution is negative for some time. It can be also found some nonnegative initial condition of sufficiently large norm at t=0 such that some component of the solution is negative at $t = 0^+$. Thus, Property (iii) is proven. □

The following result is obvious from the proof of Theorem 4.1:

**Corollary 4.2**. Theorem 4.1 (iii) is satisfied also independent of the delays for any given set of delays satisfying the constraint $0 = h_0 \leq h_1 \leq h_2 \leq ... \leq h_p = h < \infty$.

**Proof**: It follows directly since Theorem 4.1 is an independent of the delay size type result and, under the delay constraint $0 = h_0 \leq h_1 \leq h_2 \leq ... \leq h_p = h < \infty$, it has also to be fulfilled for any combination of delays satisfying the stronger constraint $0 = h_0 < h_1 < h_2 < ... < h_p = h < \infty$. □

**Corollary 4.3**. Any solution (3.8), subject to (3.9)-(3.10), to the Caputo fractional differential system (3.1) under the delay constraint $0 = h_0 \leq h_1 \leq h_2 \leq ... \leq h_p = h < \infty$ is nonnegative independent of the delays within a first interval, i.e. it satisfies $x_\alpha(t) \in R_{0+}^n$; $\forall t \in [-h, \bar{t}) \cap R_{0+}$ for some sufficiently small $\bar{t} \in \mathbf{R}_{0+}$ for any given absolutely continuous functions of initial conditions $\varphi_j : [-h, 0] \to \mathbf{R}^n_{0+}$, $j \in \overline{k-1} \cup \{0\}$ and any given piecewise continuous vector function $u: \mathbf{R}_{0+} \to \mathbf{R}^n_{0+}$ with $k=[\alpha]+1$ if $\alpha \notin \mathbf{Z}_+$ and $k = \alpha \in \mathbf{Z}_+$, $\forall t \in \mathbf{R}_{0+}$; $\forall \alpha \in \mathbf{R}_+$ if and only if $A_0 \in M\mathbf{R}^{n \times n}$, $A_i \in \mathbf{R}^{n \times n}_{0+}$ ($\forall i \in \overline{p}$) and $B \in \mathbf{R}^{n \times m}_{0+}$. In addition, and $x_\alpha : [-h, 0) \cup \mathbf{R}_{0+} \to \mathbf{R}^n_{0+}$ if, in addition, either $A_0 \geq 0$ or if $A_0$ is nilpotent or if $0 < \alpha \leq k = 1$. Furthermore, $\Psi_{\alpha j0}(t) > 0$ (with at least n entries being positive), $det \Psi_{\alpha j0}(t) > 0$ ($\forall j \in \overline{k-1} \cup \{0\}$) and $\Psi_\alpha(t) \geq 0$; $\forall t \in \mathbf{R}_{0+}$ (if $B \in \mathbf{R}^{n \times m}_+$ then $\Psi_\alpha(t) > 0$; $\forall t \in \mathbf{R}_{0+}$).

**Proof**: The solution (3.8) is identical to the unique solution (3.2) for (3.1) thus it is everywhere nonnegative under the same conditions that those of Theorem 4.1 which have been extended in Corollary 4.2. □

Note that the conditions of nonnegativity of the solution of the above theorem also imply the excitability of all the components of the state-trajectory solution; i.e., its strict positivity for some $t \in \mathbf{R}_+$ provided that $B \gg 0$ and the control $u: \mathbf{R}_{0+} \to \mathbf{R}^n_{0+}$ is admissible (i.e. piecewise continuous) and non-identically zero since $\Psi_\alpha(t) > 0$ and nonsingular for all $t \in \mathbf{R}_+$. It is now seen that the positivity conditions for the Riemann-Liouville fractional differential system (3.5) are not guaranteed in general by the above results for any given absolutely continuous functions of initial conditions $\varphi_j : [-h, 0] \to \mathbf{R}^n_{0+}$, $j \in \overline{k-1} \cup \{0\}$ and any given piecewise continuous vector function $u: \mathbf{R}_{0+} \to \mathbf{R}^n_{0+}$ with $k=[\alpha]+1$ if $\alpha \notin \mathbf{Z}_+$ and $k = \alpha \in \mathbf{Z}_+$, $\forall t \in \mathbf{R}_{0+}$; $\forall \alpha \in \mathbf{R}_+$. The following two results hold by using Corollary 3.2 and Corollary 3.4:

**Theorem 4.4**. Any solution (3.7), subject to (3.3)-(3.4), to the Riemann-Liouville fractional differential system (3.5) under the delay constraint $0 = h_0 \leq h_1 \leq h_2 \leq ... \leq h_p = h < \infty$ is everywhere nonnegative independent of the delays, i.e. it satisfies $x_\alpha : [-h, 0) \cup \mathbf{R}_{0+} \to \mathbf{R}^n_{0+}$, for any given absolutely continuous functions of initial conditions $\varphi_j : [-h, 0] \to \mathbf{R}^n_{0+}$, $j \in \overline{k-1} \cup \{0\}$ and any given piecewise continuous vector function $u: \mathbf{R}_{0+} \to \mathbf{R}^n_{0+}$ with $k=[\alpha]+1$ if $\alpha \notin \mathbf{Z}_+$ and $k = \alpha \in \mathbf{Z}_+$, $\forall t \in \mathbf{R}_{0+}$; $\forall \alpha \in \mathbf{R}_+$ if $A_0 \in M\mathbf{R}^{n \times n}$, $A_i \in \mathbf{R}^{n \times n}_{0+}$ ($\forall i \in \overline{p}$), $\left( E_{\alpha, j+1}(t) - \frac{1}{j!} I_n \right) \geq 0$; $\forall j \in \overline{k-1} \cup \{0\}$, $\forall t \in \mathbf{R}_{0+}$ and $B \in \mathbf{R}^{n \times m}_{0+}$.

The conditions $A_0 \in M\mathbf{R}^{n \times n}$, $A_i \in \mathbf{R}^{n \times n}_{0+}$ ($\forall i \in \overline{p}$) and $B \in \mathbf{R}^{n \times m}_{0+}$ are also necessary for $x_\alpha : [-h, 0) \cup \mathbf{R}_{0+} \to \mathbf{R}^n_{0+}$ for any nonnegative function of initial conditions and nonnegative controls. The condition $\left( E_{\alpha, j+1}(t) - \frac{1}{j!} I_n \right) \geq 0$; $\forall j \in \overline{k-1} \cup \{0\}$, $\forall t \in \mathbf{R}_{0+}$ is removed for initial conditions $\varphi_j : [-h, 0] \to \mathbf{R}^n_{0+}$ subject to $\varphi_j(0) = x_{j0} = 0$.

**Proof**: The proof follows in a similar way as the sufficiency part of the proof of Theorem 4.1 (iii) by inspecting the nonnegative of the solution Corollary 3.2, Eq. (3.7) for an nonnegative function of initial conditions and any nonnegative control. □

**Theorem 4.5**. Any solution (3.11), subject to (3.3)-(3.4), to the Riemann- Liouville fractional differential system (3.5) under the delay constraint $0 = h_0 \leq h_1 \leq h_2 \leq ... \leq h_p = h < \infty$ is everywhere nonnegative independent of the delays, i.e. it satisfies $x_\alpha : [-h, 0) \cup \mathbf{R}_{0+} \to \mathbf{R}_{0+}^n$, for any given absolutely continuous functions of initial conditions $\varphi_j : [-h, 0] \to \mathbf{R}_{0+}^n$, $j \in \overline{k-1} \cup \{0\}$ and any given piecewise continuous vector function $u: \mathbf{R}_{0+} \to \mathbf{R}_{0+}^n$ with $k = [\alpha]+1$ if $\alpha \notin \mathbf{Z}_+$ and $k = \alpha \in \mathbf{Z}_+$, $\forall t \in \mathbf{R}_{0+}$; $\forall \alpha \in \mathbf{R}_+$ if and only if

$$A_0 \in M\mathbf{R}^{n \times n}, \quad A_i \in \mathbf{R}_{0+}^{n \times n} \ (\forall i \in \overline{p}), \left(\Psi_{\alpha j 0}(t) - \frac{t^j}{j!} I_n\right) \geq 0; \ \forall j \in \overline{k-1} \cup \{0\}, \ \forall t \in \mathbf{R}_{0+} \text{ and } B \in \mathbf{R}_{0+}^{n \times m}.$$

The condition $\left(\Psi_{\alpha j 0}(t) - \frac{t^j}{j!} I_n\right) \geq 0; \forall j \in \overline{k-1} \cup \{0\}, \forall t \in \mathbf{R}_{0+}$ is removed for initial conditions $\varphi_j : [-h, 0] \to \mathbf{R}_{0+}^n$ subject to $\varphi_j(0) = x_{j0} = 0$.

**Proof**: The proof of sufficiency follows in a similar way as the sufficiency part of the proof of Theorem 4.1 (iii) ( see also the proof of Theorem 4.5) by inspecting the nonnegative of the solution Corollary 3.2, Eq. (3.7) for an nonnegative function of initial conditions and any nonnegative control. The proof necessity follows by contradiction by inspecting the solution (3.11) as follows:

a) Assume that $A_0 \notin M\mathbf{R}^{n \times n}$ and the solution is nonnegative for all time for any nonnegative function of initial conditions and controls. Take initial conditions $\varphi_j(t) = 0$; $\forall t \in [-h, 0]$, $\forall j \in \overline{k-1}$; $\varphi_0(t) = 0$; $\forall t \in [-h, 0)$, $\varphi_0(0) = x_{00}(0) \neq 0$ and $u \equiv 0$ on $\mathbf{R}_{0+}$. Then (3.11) becomes:

$$x_\alpha(t) = (\Psi_{\alpha 00}(t) - I_n) x_{00} = (E_{\alpha 00}(t) - I_n) x_{00} \text{ for } t \in [0, h_1]$$

since $\Psi_{\alpha 00}(t) = \Phi_{\alpha 00}(t)$ for $t \in [0, h_1]$. Since $A_0 \notin M\mathbf{R}^{n \times n}$, $\exists t \in [0, h_1]$ and $(i = i(t), m(t) \neq i) \in \overline{n}$ such that $(\Phi_{\alpha 00}(t))_{im} < 0$. Otherwise, if $A_0 \notin M\mathbf{R}^{n \times n}$ and $\Phi_{\alpha 00}(t) > 0$; $\forall t \in [0, h_1]$, it would follow from (4.3) that $\Phi_{\alpha 00}(t) > 0$; $\forall t \in \mathbf{R}_{0+}$ since

$$e^{A_0 t} = e^{\chi A_0 h_1 + A_0 \delta} = e^{\chi A_0 h_1} e^{A_0 \delta} = \left(e^{A_0 h_1}\right)^\chi e^{A_0 \delta} > 0$$

from the semigroup property of $\left(e^{A_0 t}, t \in \mathbf{R}_{0+}\right)$ with $\chi = \chi(t, h_1) = [t/h_1]$ and $(0, h_1] \ni \delta = \delta(t, h_1) = t - \chi h_1$ what implies $\Phi_{\alpha 00}(t) > 0$; $\forall t \in \mathbf{R}_{0+}$ from (4.3). Thus, $A_0 \in M\mathbf{R}^{n \times n}$ which contradicts $A_0 \notin M\mathbf{R}^{n \times n}$. It has been proven that $A_0 \notin M\mathbf{R}^{n \times n} \Rightarrow e_i^T \Phi_{\alpha 00}(t) e_m < 0$; $\forall t \in (0, h_1]$ for some $(i = i(t), m(t) \neq i) \in \overline{n}$. Now, take $x_{00j} = \delta_{jm} (\forall j \in \overline{n})$ where $\delta_{jm}$ denotes the Kronecker delta. Then,

$$x_{\alpha i}(t) = e_i^T\left(\Psi_{\alpha 00}(t) - \frac{t^j}{j!}I_n\right)x_{00} = e_i^T\left(\Phi_{\alpha 00}(t) - \frac{t^j}{j!}I_n\right)e_m = e_i^T\Phi_{\alpha 00}(t)e_m < 0$$

As a result, $A_0 \in MR^{n \times n}$ is a necessary condition for the solution to be nonnegative for all time irrespective of the delay sizes.

b) Assume that the solution is nonnegative for all time for any nonnegative function of initial conditions and controls Assume that $e_i^T A_\ell e_j < 0$ and $h_\ell \neq h_i; \forall i(\neq \ell) \in \bar{p}$ for some $i,j \in \bar{n}$, $\ell \in \bar{p}$. Take initial conditions $x_{j0} = \varphi_j(0) = 0$; $\forall t \in [-h, 0]$; $\forall j \in \overline{k-1} \cup \{0\}$, $\varphi_j \equiv 0$; $\forall j \in \overline{k-1}$ and $u \equiv 0$. One gets from (3.2)

$$x_{\alpha i}(t) = \int_0^{h'} e_i^T \Phi_\alpha(t-\tau)\left(\sum_{j(\neq\ell)\in\bar{p}} A_j\right)\varphi_0(\tau - h')d\tau + \int_0^{h_\ell} e_i^T \Phi_\alpha(t-\tau)A_\ell \varphi_0(\tau - h_\ell)d\tau; \forall t \in [0, h_i]$$

for the case $h' = h_i$; $\forall i(\neq \ell \in \bar{p})$. Now, if $h = h_\ell > h'$, Take a further specification of initial conditions as follows: $\varphi_0(t) = 0$; $\forall t \in [0, h']$, and $\varphi_0(\tau) = (k_1, ..., k_n)^T >> 0$; $\forall t \in (h', h_\ell]$ then

$$x_{\alpha i}(t) = \int_{h'}^{h_\ell} e_i^T \Phi_\alpha(t-\tau)A_\ell \varphi_0(\tau - h_\ell)d\tau = \sum_{r=1}^n \sum_{m=1}^n \left(\int_{h'}^{h_\ell} \Phi_{\alpha ir}(t-\tau)A_{\ell rm}d\tau\right)k_m$$

$$= \sum_{j=1}^n k_j \left(\int_{h'}^{h_\ell} \Phi_{\alpha i}^T(t-\tau)d\tau\right)A_{\ell j} = \sum_{(m\neq j)\in\bar{n}} k_m \left(\sum_{r=1}^m \int_{h'}^{h_\ell} \Phi_{\alpha ir}(t-\tau)d\tau\right)A_{\ell rm} + k_j \left(\sum_{r=1}^m \int_{h'}^{h_\ell} \Phi_{\alpha ir}(t-\tau)d\tau\right)A_{\ell rj}$$

; $\forall t \in [0, h_i]$. As a result, $A_i \geq 0 (\forall i \in \bar{p})$ is a necessary condition for the solution to be nonnegative for all time irrespective of the delay sizes.

c) Assume that the solution is nonnegative for all time for any nonnegative function of initial conditions and controls and $B \geq 0$ is not fulfilled so that it exists at least an entry $B_{\ell j} < 0$ of B. Then, one has under identically zero initial conditions the following unique solution:

$$x_{\alpha i}(t) = \int_0^t e_i^T \Psi_\alpha(t-\tau)Bu(\tau)d\tau = \sum_{i=1}^m \int_0^t \Psi_{\alpha i}^T(t-\tau)B_i u_i(\tau)d\tau = \sum_{i=1}^m \sum_{\ell=1}^n \int_0^t \Psi_{\alpha i\ell}(t-\tau)B_{\ell i}u_i(\tau)d\tau$$

$$= \sum_{(i\neq j)\in\bar{m}} \sum_{\ell=1}^n \int_0^t \Psi_{\alpha i\ell}(t-\tau)B_{\ell i}u_i(\tau)d\tau - \left(\sum_{\ell=1}^n \int_0^t \Psi_{\alpha j\ell}(t-\tau)|B_{\ell j}|d\tau\right)k_{uj} < 0; t \in R_{0+}$$

provided that

$$k_{uj} > \frac{\sum_{(i\neq j)\in\bar{m}} \sum_{\ell=1}^n \int_0^t \Psi_{\alpha i\ell}(t-\tau)B_{\ell i}u_i(\tau)d\tau}{\left(\sum_{\ell=1}^n \int_0^t \Psi_{\alpha j\ell}(t-\tau)|B_{\ell j}|d\tau\right)}$$

by assuming that $B \geq 0$ fails because $B_{\ell j} < 0$ for some $(\ell, j) \in \bar{n} \times \bar{m}$ and a constant control component $u_j \equiv k_{uj} > 0$ is injected on the time interval $[0, t]$ for some arbitrary $t \in R_+$ for the remaining control components being chosen t be nonnegative for all time. This contradicts that the solution is nonnegative for all time if the condition $B \geq 0$ fails. □

*Remark 4.6.* Note that Theorem 4.1 can be extended as a necessary condition for $t \in [0, h_1]$ since $\Psi_{\alpha j 0}(t) = \Phi_{\alpha j 0}(t)$ for $t \in [0, h_1]$; $\forall j \in \overline{k-1} \cup \{0\}, \forall t \in R_{0+}$. □

*Remark 4.7.* Note by simple calculation that $\left(A_0 \in MR^{n \times n} \wedge A_i \geq 0 (i \in \overline{p})\right) \Rightarrow \left(\sum_{i=0}^{p} A_i\right) \in MR^{n \times n}$. This is a necessary and sufficient condition for the nonnegativity of the solutions of the Caputo fractional differential system (3.1) of arbitrary order $\alpha \in R_+$ under arbitrary nonnegative controls and initial conditions in the absence of delays; i.e. for $h_i = 0$; $\forall i \in \overline{\omega} \cup \{0\}$ and any $\omega \in Z_+$. □

*Remark 4.8.* The given conditions to guarantee that the solution is everywhere nonnegative under any given arbitrary nonnegative initial conditions and nonnegative controls are of independent of the sizes of the delays type; i.e. for any given set of p delays. However, the conditions are weakened for particular situations involving repeated delays as follows. Note from Theorem 4.5 that the various given conditions $A_i \geq 0$ of necessary type to guarantee the nonnegativity of the solution under any admissible nonnegative controls and nonnegative initial conditions are weakened to $\left(\sum_{\ell=1}^{\nu_i} A_{\left(\sum_{j=0}^{i-1} \nu_j + \ell\right)}\right)$ if there is some repeated delay $h_i$ of multiplicity $\nu_i \geq 2$ (i.e. the number of distinct delays is $0 \leq q < p$). Also, if $h_0 = 0$ is repeated with multiplicity $\nu_0 \geq 2$ then the condition $A_0 \in MR^{n \times n}$ for $\nu_0 = 1$ is replaced by $\left(\sum_{\ell=0}^{\nu_0 - 1} A_\ell\right) \in MR^{n \times n}$. □

*Remark 4.9.* Note that there is a duality of all the given results of sufficiency type or necessary and sufficiency type in the sense that the solutions are guaranteed to be nonpositive for all time under similar conditions for the cases when all components of the controls and initial conditions are nonpositive for all time. □

## 5. Asymptotic behavior of unforced solutions for $\alpha \in R_+$

The asymptotic behaviour and the stability properties of the Caputo fractional differential system (3.1) can be investigated via the extension of the subsequent formulas for $\alpha \in R_+$, (see (1.8.27)-(1.8.29), [1]):

1) If $0 < \alpha < 2$ then for $|z| \to \infty$ and some $\mu \in R$ satisfying $\mu < \pi \min(1, \alpha)$:

$$E_{\alpha \beta}(z) = \frac{1}{\alpha} z^{(1-\beta)/\alpha} e^{\left(z^{1/\alpha}\right)} - \sum_{j=1}^{N} \frac{1}{\Gamma(\beta - \alpha j)} \frac{1}{z^j} + O\left(\frac{1}{z^{N+1}}\right) \quad (5.1)$$

with $|arg \, z| \leq \mu < \pi \min(1, \alpha)$, any $N \in Z_+$, and

$$E_{\alpha \beta}(z) = -\sum_{j=1}^{N} \frac{1}{\Gamma(\beta - \alpha j)} \frac{1}{z^j} + O\left(\frac{1}{z^{N+1}}\right) \quad (5.2)$$

with $\pi \geq |arg\, z| \geq \mu < \pi = \pi\, min(1, \alpha)$, any $N \in \mathbf{Z}_+$

2) If $\alpha \geq 2$ then for $|z| \to \infty$

$$E_{\alpha\beta}(z) = \frac{1}{\alpha} \sum_{j \in Q} \left( z^{1/\alpha} e^{\frac{2j\pi i}{\alpha}} \right)^{1-\beta} e^{\left( e^{\frac{2j\pi i}{\alpha}} z^{1/\alpha} \right)} - \sum_{j=1}^{N} \frac{1}{\Gamma(\beta - \alpha j)} \frac{1}{z^j} + O\left( \frac{1}{z^{N+1}} \right) \quad (5.3)$$

for any $N \in \mathbf{Z}_+$ with $\beta \in \overline{k-1} \cup \{0\} \cup \{\alpha\}$, $|arg\, z| \leq \frac{\alpha\pi}{2}$, $Q := \left\{ n \in \mathbf{Z} : |arg\, z + 2\pi n| \leq \frac{\alpha\pi}{2} \right\}$ and

$i = \sqrt{-1}$ being the complex imaginary unit. The above formulas are extendable to the Mittag-Leffler matrix functions $E_{\alpha, j+1}(A_0 t^\alpha) := \sum_{\ell=0}^{\infty} \frac{(A_0 t^\alpha)^\ell}{\Gamma(\alpha\ell + j + 1)}$; $\forall j \in \overline{k-1} \cup \{0\}$, respectively, $E_{\alpha\alpha}(A_0 t^\alpha)$ by identifying $z \to A_0 t^\alpha$, $z^{\frac{1}{\alpha}} \to (A_0)^{\frac{1}{\alpha}} t$ (if $(A_0)^{\frac{1}{\alpha}}$ exists) and $z^{-1} \to A_0^{-1} t^{-\alpha}$ (if $A_0$ is non-singular),

$\beta \to j+1$, respectively, $\beta \to \alpha$. Irrespective of the existence of $(A_0)^{\frac{1}{\alpha}}$ and of $A_0 \neq 0$ being singular or nonsingular, it is possible to identify $z^{-1} \to \|A_0\|^{-1} t^{-\alpha}$ and $z \to \|A_0\| t^\alpha$ and to use

$$\|E_{\alpha, j+1}(A_0 t^\alpha)\| \leq E_{\alpha, j+1}(\|A_0\| t^\alpha) = \sum_{\ell=0}^{\infty} \frac{(\|A_0\| t^\alpha)^\ell}{\Gamma(\alpha\ell + j + 1)}; \quad \forall j \in \overline{k-1} \cup \{0\} \quad (5.4)$$

$$\|E_{\alpha\alpha}(A_0 t^\alpha)\| \leq E_{\alpha\alpha}(\|A_0\| t^\alpha) = \sum_{\ell=0}^{\infty} \frac{(\|A_0\| t^\alpha)^\ell}{\Gamma(\alpha\ell + \alpha)} \quad (5.5)$$

The method may be used to calculate an asymptotic estimate of the solution (3.2) if $A_0$ is non-singular (or an upper-bounding function for any nonzero $A_0$) of the Caputo fractional differential system (3.1), via (3.3)-(3.4), or, equivalently (3.8), via (3.9)-(3.10) and (3.3)-(3.4). The estimations may be extended with minor modification to the Riemann-Liouville fractional differential system (3.5). Note that if all the complex eigenvalues of $A_0$ appear by conjugate pairs $A_0$ then $A_0 = T^{-1} J_0 A_0$ where $J_0$ is its real canonical form.

A) Assume that $\alpha \in \mathbf{R}_+$, $A_0$ is real non-singular and $(A_0)^{\frac{1}{\alpha}}$ exists; i.e. $\exists M$ such that $M^\alpha = A_0$ and $A_i (i \in \overline{p})$ is real. Then, one gets from (5.1)-(5.3):

$$E_{\alpha, j+1}(A_0 t^\alpha) = \frac{1}{\alpha} \left( A_0^{-j} \right)^{\frac{1}{\alpha}} t^{-j} e^{\left( A_0^{\frac{1}{\alpha}} \right) t} - \sum_{\ell=1}^{N} \frac{1}{\Gamma((1-\alpha)\ell+1)} \left( A_0^{-\ell} \right) t^{-\ell\alpha} + O\left( A_0^{-(N+1)} t^{-(N+1)\alpha} \right)$$

$$; \forall j \in \overline{k-1} \cup \{0\} \quad (5.6)$$

$$E_{\alpha\alpha}(A_0 t^\alpha) = \frac{1}{\alpha} A_0^{(1-\alpha)/\alpha} t^{1-\alpha} e^{\left( A_0^{\frac{1}{\alpha}} \right) t} - \sum_{\ell=1}^{N} \frac{1}{\Gamma((1-\ell)\alpha)} \left( A_0^{-\ell} \right) t^{-\ell\alpha} + O\left( A_0^{-(N+1)} t^{-(N+1)\alpha} \right)$$

$$(5.7)$$

as $t \to \infty$ if $0 < \alpha < 2$, for any $N \in \mathbf{Z}_+$, and

$$E_{\alpha, j+1}\left(A_0 t^\alpha\right) = \frac{1}{\alpha} \sum_{\ell \in Q} \left(A_0^{1/\alpha} t \, e^{\frac{2\ell \pi i}{\alpha}}\right)^{-j} e^{\left(e^{\frac{2\ell \pi i}{\alpha}}\left(A_0^{1/\alpha}\right) t\right)}$$

$$- \sum_{\ell=1}^{N} \frac{1}{\Gamma((1-\alpha)\ell+1)} \left(A_0^{-\ell}\right) t^{-\ell \alpha} + O\left(A_0^{-(N+1)} t^{-(N+1)\alpha}\right) ; \forall j \in \overline{k-1} \cup \{0\} \quad (5.8)$$

$$E_{\alpha \alpha}\left(A_0 t^\alpha\right) = \frac{1}{\alpha} \sum_{\ell \in Q} \left(A_0^{1/\alpha} t \, e^{\frac{2\ell \pi i}{\alpha}}\right)^{1-\alpha} e^{\left(e^{\frac{2\ell \pi i}{\alpha}}\left(A_0^{1/\alpha}\right) t\right)}$$

$$- \sum_{\ell=1}^{N} \frac{1}{\Gamma((1-\ell)\alpha)} \left(A_0^{-\ell}\right) t^{-\ell \alpha} + O\left(A_0^{-(N+1)} t^{-(N+1)\alpha}\right) \quad (5.9)$$

as $t \to \infty$ if $\alpha \geq 2$, for any $N \in \mathbf{Z}_+$, with $Q := \left\{ n \in \mathbf{Z} : |n| \leq \frac{\alpha}{4} \right\}$

B) Assume that $\alpha \in \mathbf{R}_+$ and $A_i$ $\left(i \in \overline{p} \cup \{0\}\right)$ is real, one obtains from (5.1)-(5.2):

$$\left|E_{\alpha \beta}(z)\right| \leq \overline{E}_{\alpha \beta}(z) = \frac{1}{\alpha} |z|^{(1-\beta)/\alpha} \left| e^{z^{1/\alpha}} + \left|\sum_{\ell=1}^{N} \frac{1}{\Gamma(\beta - \alpha \ell)} \frac{1}{z^j} + O\left(\frac{1}{z^{N+1}}\right)\right| \right| \quad (5.10)$$

for $0 < \alpha < 2$, $\beta \in \overline{k-1} \cup \{0\}$, and

$$\left|E_{\alpha \beta}(z)\right| \leq \overline{E}_{\alpha \beta}(z) = \frac{1}{\alpha} \sum_{j \in Q} \left(|z|^{(1-\beta)/\alpha}\right) \left| e^{z^{1/\alpha}} + \left|\sum_{\ell=1}^{N} \frac{1}{\Gamma(\beta - \alpha \ell)} \frac{1}{z^j} + O\left(\frac{1}{z^{N+1}}\right)\right| \right| \quad (5.11)$$

for $\alpha \geq 2$, $\beta \in \overline{k-1} \cup \{0, \alpha\}$. Thus, on gets from (5.10)

$$\left\|E_{\alpha, j+1}\left(A_0 t^\alpha\right)\right\| \leq \frac{1}{\alpha} \left(\|A_0\|^{-j}\right)^{\frac{1}{\alpha}} t^{-j} \left\|e^{A_0 t}\right\|^{\frac{1}{\alpha}} - \sum_{\ell=1}^{N} \frac{1}{\Gamma((1-\alpha)\ell+1)} \left(\|A_0\|^{-\ell}\right) t^{-\ell \alpha} + O\left(\|A_0\|^{-(N+1)} t^{-(N+1)\alpha}\right)$$

$$; \forall j \in \overline{k-1} \cup \{0\} \quad (5.12)$$

$$\left\|E_{\alpha \alpha}\left(A_0 t^\alpha\right)\right\| \leq \frac{1}{\alpha} \|A_0\|^{(1-\alpha)/\alpha} t^{1-\alpha} \left\|e^{A_0 t}\right\|^{\frac{1}{\alpha}} - \sum_{\ell=1}^{N} \frac{1}{\Gamma((1-\ell)\alpha)} \left(\|A_0\|^{-\ell}\right) t^{-\ell \alpha} + O\left(\|A_0\|^{-(N+1)} t^{-(N+1)\alpha}\right)$$

$$(5.13)$$

as $t \to \infty$, for any $N \in \mathbf{Z}_+$, if $0 < \alpha < 2$, and one gets from (5.11)

$$\left\|E_{\alpha, j+1}\left(A_0 t^\alpha\right)\right\| \leq \frac{1}{\alpha} \sum_{\ell \in Q} \left(\|A_0\|^{1/\alpha} t\right)^{-j} \left\|e^{A_0 t}\right\|^{\frac{1}{\alpha}}$$

$$- \sum_{\ell=1}^{N} \frac{1}{\Gamma((1-\alpha)\ell+1)} \left(\|A_0\|^{-\ell}\right) t^{-\ell \alpha} + O\left(\|A_0\|^{-(N+1)} t^{-(N+1)\alpha}\right) ; \forall j \in \overline{k-1} \cup \{0\} \quad (5.14)$$

$$\left\|E_{\alpha \alpha}\left(A_0 t^\alpha\right)\right\| \leq \frac{1}{\alpha} \sum_{\ell \in Q} \left(\|A_0\|^{1/\alpha} t\right)^{1-\alpha} \left\|e^{A_0 t}\right\|^{\frac{1}{\alpha}}$$

$$-\sum_{\ell=1}^{N}\frac{1}{\Gamma((1-\ell)\alpha)}\left(\|A_0\|^{-\ell}\right)t^{-\ell\alpha}+O\left(\|A_0\|^{-(N+1)}t^{-(N+1)\alpha}\right) \tag{5.15}$$

as $t\to\infty$ if $\alpha\geq 2$, for any $N\in\mathbf{Z}_+$, with $Q:=\left\{n\in\mathbf{Z}:|n|\leq\dfrac{\alpha}{4}\right\}$. The formula (3.8) for the solution is more useful than its equivalent expression (3.2) to investigate the asymptotic properties of the Caputo fractional differential system. Therefore, we obtain now either explicit or upper-bounding asymptotic expressions for (3.9)-(3.10) by using (5.6) to (5.13) as follows:

1) Assume that $\alpha\in\mathbf{R}_+$, $A_0$ is real non-singular, $(A_0)^{\frac{1}{\alpha}}$ exists and $A_i$ $(i\in\overline{p})$ are also real. Then, one gets from (5.6)-(5-9) into (3.9)-(3.10):

$$\Psi_{\alpha j 0}(t)=\frac{1}{\alpha}\left(A_0^{-j}\right)^{\frac{1}{\alpha}}e^{\left(A_0^{1/\alpha}\right)t}-\sum_{\ell=1}^{N}\frac{1}{\Gamma(j+1-\alpha\ell)}\left(A_0^{-\ell}\right)t^{j-\ell\alpha}+O\left(A_0^{-(N+1)}t^{j-(N+1)\alpha}\right)$$

$$+\sum_{i=1}^{p}\int_0^t\left(\frac{1}{\alpha}A_0^{(1-\alpha)/\alpha}e^{\left(A_0^{1/\alpha}\right)\tau}-\sum_{\ell=1}^{N}\frac{1}{\Gamma((1-\ell)\alpha)}\left(A_0^{-\ell}\right)\tau^{(1-\ell)\alpha-1}+O\left(A_0^{-(N+1)}\tau^{-N\alpha-1}\right)\right)A_i\Psi_{\alpha j 0}(t-\tau-h_i)d\tau$$

(5.16)

$$\Psi_\alpha(t)=\frac{1}{\alpha}A_0^{(1-\alpha)/\alpha}e^{\left(A_0^{1/\alpha}\right)t}-\sum_{\ell=1}^{N}\frac{1}{\Gamma((1-\ell)\alpha)}\left(A_0^{-\ell}\right)t^{(1-\ell)\alpha}+O\left(A_0^{-(N+1)}t^{-N\alpha}\right)$$

$$+\sum_{i=1}^{p}\int_0^t\left(\frac{1}{\alpha}A_0^{(1-\alpha)/\alpha}e^{\left(A_0^{1/\alpha}\right)\tau}-\sum_{\ell=1}^{N}\frac{1}{\Gamma((1-\ell)\alpha)}\left(A_0^{-\ell}\right)\tau^{(1-\ell)\alpha-1}+O\left(A_0^{-(N+1)}\tau^{-N\alpha-1}\right)\right)A_i\Psi_\alpha(t-\tau-h_i)d\tau$$

(5.17)

$\forall j\in\overline{k-1}\cup\{0\}$ as $t\to\infty$ if $0<\alpha<2$, for any $N\in\mathbf{Z}_+$, and

$$\Psi_{\alpha j 0}(t)=\frac{1}{\alpha}\sum_{\ell\in Q}\left(A_0^{1/\alpha}e^{\frac{2\ell\pi i}{\alpha}}\right)^{-j}e^{\left(e^{\frac{2\ell\pi i}{\alpha}}\left(A_0^{1/\alpha}\right)t\right)}-\sum_{\ell=1}^{N}\frac{1}{\Gamma(j+1-\alpha\ell)}\left(A_0^{-\ell}\right)t^{j-\ell\alpha}+O\left(A_0^{-(N+1)}t^{j-(N+1)\alpha}\right)$$

$$+\sum_{i=1}^{p}\int_0^t\left(\frac{1}{\alpha}\sum_{\ell\in Q}\left(A_0^{(1-\alpha)/\alpha}e^{\frac{2\ell\pi i}{\alpha}}\right)^{-j}e^{\left(e^{\frac{2\ell\pi i}{\alpha}}\left(A_0^{1/\alpha}\right)\tau\right)}-\sum_{\ell=1}^{N}\frac{1}{\Gamma((1-\ell)\alpha)}\left(A_0^{-\ell}\right)\tau^{(1-\ell)\alpha-1}+O\left(A_0^{-(N+1)}\tau^{-N\alpha-1}\right)\right)A_i\Psi_{\alpha j 0}(t-\tau-h_i)d\tau$$

(5.18)

$$\Psi_\alpha(t)=\frac{1}{\alpha}\sum_{\ell\in Q}\left(A_0^{(1-\alpha)/\alpha}e^{\frac{2\ell\pi(1-\alpha)i}{\alpha}}\right)e^{\left(e^{\frac{2\ell\pi i}{\alpha}}\left(A_0^{1/\alpha}\right)t\right)}-\sum_{\ell=1}^{N}\frac{1}{\Gamma((1-\alpha)\ell+1)}\left(A_0^{-\ell}\right)t^{j-\ell\alpha}+O\left(A_0^{-(N+1)}t^{j-(N+1)\alpha}\right)$$

$$+\sum_{i=1}^{p}\int_0^t\left(\frac{1}{\alpha}\sum_{\ell\in Q}\left(A_0^{(1-\alpha)/\alpha}e^{\frac{2\ell\pi(1-\alpha)i}{\alpha}}\right)e^{\left(e^{\frac{2\ell\pi i}{\alpha}}\left(A_0^{1/\alpha}\right)\tau\right)}-\sum_{\ell=1}^{N}\frac{1}{\Gamma((1-\alpha)\ell)}\left(A_0^{-\ell}\right)\tau^{(1-\ell)\alpha-1}+O\left(A_0^{-(N+1)}\tau^{-N\alpha-1}\right)\right)A_i\Psi_{\alpha j 0}(t-\tau-h_i)d\tau$$

(5.19)

$\forall j \in \overline{k-1} \cup \{0\}$ as $t \to \infty$ if $\alpha \geq 2$, for any $N \in \mathbf{Z}_+$.

2) Assume that $\alpha \in \mathbf{R}_+$ and $A_i \, (i \in \overline{p} \cup \{0\})$ are real. Then,

$$\|\Psi_{\alpha j 0}(t)\| \leq \frac{1}{\alpha} \left(\|A_0\|^{-j}\right)^{\frac{1}{\alpha}} \|e^{A_0 t}\|^{\frac{1}{\alpha}} + \sum_{\ell=1}^{N} \frac{1}{|\Gamma(j+1-\alpha \ell)|} \|A_0\|^{-\ell} t^{j-\ell \alpha} + O\left(\|A_0\|^{-(N+1)} t^{j-(N+1)\alpha}\right)$$

$$+ \sum_{i=1}^{p} \int_0^t \left( \frac{1}{\alpha} \|A_0\|^{(1-\alpha)/\alpha} \|e^{A_0 \tau}\|^{\frac{1}{\alpha}} + \sum_{\ell=1}^{N} \frac{1}{|\Gamma((1-\ell)\alpha)|} \left(\|A_0\|^{-\ell}\right) \tau^{(1-\ell)\alpha - 1} + O\left(\|A_0\|^{-(N+1)} \tau^{-N\alpha - 1}\right) \right) \|A_i\| \|\Psi_{\alpha j 0}(t-\tau-h_i)\| d\tau$$

(5.20)

$$\|\Psi_\alpha(t)\| \leq \frac{1}{\alpha} \|A_0\|^{(1-\alpha)/\alpha} \|e^{A_0 t}\|^{\frac{1}{\alpha}} t^\alpha + \sum_{\ell=1}^{N} \frac{1}{|\Gamma((1-\ell)\alpha)|} \left(\|A_0\|^{-\ell}\right) t^{(1-\ell)\alpha} + O\left(\|A_0\|^{-(N+1)} t^{-N\alpha}\right)$$

$$+ \sum_{i=1}^{p} \int_0^t \left( \frac{1}{\alpha} \|A_0\|^{(1-\alpha)/\alpha} \tau^\alpha \|e^{A_0 \tau}\|^{\frac{1}{\alpha}} + \sum_{\ell=1}^{N} \frac{1}{|\Gamma((1-\ell)\alpha)|} \|A_0\|^{-\ell} \tau^{(1-\ell)\alpha - 1} + O\left(\|A_0\|^{-(N+1)} \tau^{-N\alpha - 1}\right) \right) \|A_i\| \|\Psi_\alpha(t-\tau-h_i)\| d\tau$$

(5.21)

$\forall j \in \overline{k-1} \cup \{0\}$ as $t \to \infty$ if $0 < \alpha < 2$, for any $N \in \mathbf{Z}_+$, and

$$\|\Psi_{\alpha j 0}(t)\| \leq \frac{1}{\alpha} \left(\|A_0\|^{1/\alpha}\right)^{-j} \|e^{A_0 t}\|^{\frac{1}{\alpha}} + \sum_{\ell=1}^{N} \frac{1}{|\Gamma(j+1-\alpha \ell)|} \|A_0\|^{-\ell} t^{j-\ell \alpha} + O\left(\|A_0\|^{-(N+1)} t^{j-(N+1)\alpha}\right)$$

$$+ \sum_{i=1}^{p} \int_0^t \left( \frac{1}{\alpha} \|A_0\|^{(1-\alpha)/\alpha} \|e^{A_0 \tau}\|^{\frac{1}{\alpha}} + \sum_{\ell=1}^{N} \frac{1}{|\Gamma((1-\ell)\alpha)|} \|A_0\|^{-\ell} \tau^{(1-\ell)\alpha - 1} + O\left(\|A_0\|^{-(N+1)} \tau^{-N\alpha - 1}\right) \right) \|A_i\| \|\Psi_{\alpha j 0}(t-\tau-h_i)\| d\tau$$

(5.22)

$$\|\Psi_\alpha(t)\| \leq \frac{1}{\alpha} \|A_0\|^{(1-\alpha)/\alpha} \|e^{A_0 t}\|^{\frac{1}{\alpha}} + \sum_{\ell=1}^{N} \frac{1}{|\Gamma((1-\ell)\alpha)|} \|A_0\|^{-\ell} t^{j-\ell \alpha} + O\left(\|A_0\|^{-(N+1)} t^{j-(N+1)\alpha}\right)$$

$$+ \sum_{i=1}^{p} \int_0^t \left( \frac{1}{\alpha} \|A_0\|^{(1-\alpha)/\alpha} \|e^{A_0 \tau}\|^{\frac{1}{\alpha}} + \sum_{\ell=1}^{N} \frac{1}{|\Gamma((1-\ell)\alpha)|} \|A_0\|^{-\ell} \tau^{(1-\ell)\alpha - 1} + O\left(\|A_0\|^{-(N+1)} \tau^{-N\alpha - 1}\right) \right) \|A_i\| \|\Psi_{\alpha j 0}(t-\tau-h_i)\| d\tau$$

(5.23)

$\forall j \in \overline{k-1} \cup \{0\}$ as $t \to \infty$ if $\alpha \geq 2$, for any $N \in \mathbf{Z}_+$.

For further discussion, note that it exists a set of linearly independent continuously differential real functions $\{\alpha_i : \mathbf{R}_{0+} \to \mathbf{R}, i \in \overline{\nu - 1} \cup \{0\}\}$, where $\nu$ is the degree of the minimal polynomial of any square real matrix $A_0$ such that:

$$e^{A_0 t} = \sum_{i=0}^{\nu-1} \alpha_i(t) A_0^i = \sum_{j=0}^{\nu} \sum_{i=0}^{\nu_i - 1} k_{ij} t^i e^{\lambda_j t} \,; \quad \forall t \in \mathbf{R}_{0+} \quad (5.24)$$

(see, for instance, [4-5]), where $k_{ij} \in \mathbf{R}$; $i \in \overline{\nu-1} \cup \{0\}$, $j \in \overline{\nu} \cup \{0\}$, $\sigma(M) := \{\lambda_i \in \mathbf{C} : det(\lambda_i I_n - A_0) = 0\}$ is the spectrum of $A_0$ defined by the set of eigenvalues $\lambda_i$ of M of respective index $\nu_i$ (i.e. the multiplicity of $\lambda_i$ in the minimal polynomial of $A_0$) and algebraic multiplicity $\mu_i$ (i.e. the multiplicity of $\lambda_i$ in the characteristic polynomial of $A_0$) so that $n = \sum_{i=1}^{n} n_i \geq \nu = \sum_{i=1}^{n} \nu_i$ with n being the order of $A_0$ with $\nu$ being the degree of its minimal polynomial.

The subsequent stability result is based on the above formulas:

**Theorem 5.1**. The following properties hold:

(i) If $k = \alpha = 1$ (the particular standard bon- fractional case) then (3.1) is globally Lyapunov stable independent of the delays if

$$\left\| \left( \frac{1}{\beta_1} A_1, \frac{1}{\beta_2} A_2, \cdots, \frac{1}{\beta_p} A_p \right) \right\|_2 \leq -\mu_2(A_0) \qquad (5.25)$$

requiring for the $\ell_2$- matrix measure of $A_0$ to fulfil $\mu_2(A_0) := \frac{1}{2} \lambda_{max}\left(A_0 + A_0^T\right) \leq 0$, for some $\beta_i \in \mathbf{R}_+$ ($i \in \bar{p}$) subject to $\sum_{i=1}^{p} \beta_i^2 = 1$, [6]. Also,

$$\Psi_{100}(t) = e^{A_0 t} - A_0^{-1} t^{-1} + O\left(A_0^{-1} t^{-1}\right) + \sum_{i=1}^{p} \int_0^t \left( e^{A_0 \tau} - A_0^{-1} \tau^{-1} + O\left(A_0^{-1} \tau^{-1}\right) \right) A_i \Psi_{100}(t - \tau - h_i) d\tau$$

as $t \to \infty$ (5.26)

is bounded provided that $A_0$ is non- singular with $e^{A_0 t}$ being of the form (5.24) if (5.25) holds and then the unforced solution:

$$x_\alpha(t) = \sum_{j=0}^{k-1} \left( \Psi_{\alpha j 0}(t) x_{j0} + \sum_{i=1}^{p} \int_0^{h_i} \Psi_{\alpha j 0}(t - \tau) \varphi_j(\tau - h_i) d\tau \right) \qquad (5.27)$$

Is bounded for all time. Furthermore,

$$\|\Psi_{100}(t)\| \leq \left\| e^{A_0 t} \right\| + \|A_0\|^{-1} t^{-1} + O\left(\|A_0\|^{-2} t^{-2}\right)$$

$$+ \sum_{i=1}^{p} \int_0^t \left( \left\| e^{A_0 \tau} \right\| + \|A_0\|^{-1} \tau^{-1} + O\left(\|A_0\|^{-1} \tau^{-1}\right) \right) \|A_i\| \|\Psi_{100}(t - \tau - h_i)\| d\tau \qquad \text{as } t \to \infty \qquad (5.28)$$

if (5.25) holds irrespective of $A_0$ being singular or non-singular. If, in addition, $\mu_2(A_0) < 0$ and (5.25) holds with strict inequality then (3.1) is globally asymptotically Lyapunov stable independent of the delays and

$$\Psi_{100}(t) \to \sum_{i=1}^{p} \int_0^t \left( e^{A_0 \tau} - A_0^{-1} \tau^{-1} + O\left(A_0^{-1} \tau^{-1}\right) \right) A_i \Psi_{10}(t - \tau - h_i) d\tau \to 0 \qquad \text{as } t \to \infty \qquad (5.29)$$

(ii) If $k = 1$ and $\alpha \in (0, 1]$ the inequality (5.25) is strict then (3.1) is globally Lyapunov stable independent of the delays if $\mu_2\left(A_0^{1/\alpha}\right) \leq 0$ and

$$\left\| \left( \frac{1}{\beta_1} A_1, \frac{1}{\beta_2} A_2, \cdots, \frac{1}{\beta_p} A_p \right) \right\|_2 \leq -\mu_2\left(A_0^{1/\alpha}\right) \qquad (5.30)$$

provided that $A_0$ is non- singular and $A_0^{1/\alpha}$ exists. Also, then (3.1) is globally asymptotically Lyapunov stable independent of the delays if, in addition, $\mu_2\left(A_0^{1/\alpha}\right) < 0$ and

$$\left\|\left(\frac{1}{\beta_1}A_1, \frac{1}{\beta_2}A_2, \cdots, \frac{1}{\beta_p}A_p\right)\right\|_2 < \left|\mu_2\left(A_0^{1/\alpha}\right)\right| \tag{5.31}$$

$$\Psi_{\alpha 00}(t) = \frac{1}{\alpha} e^{\left(A_0^{1/\alpha}\right)t} - \frac{1}{\Gamma(1-\alpha)} A_0^{-1} t^{-\alpha} + O\left(A_0^{-2} t^{-2\alpha}\right)$$

$$+ \sum_{i=1}^{p} \int_0^t \left[\frac{1}{\alpha} A_0^{(1-\alpha)/\alpha} e^{\left(A_0^{1/\alpha}\right)\tau} - A_0^{-1}\tau^{-1} + O\left(A_0^{-(N+1)}\tau^{-\alpha-1}\right)\right] A_i \Psi_{\alpha 00}(t-\tau-h_i) d\tau \to 0$$

$$\text{as } t \to \infty \tag{5.32}$$

If either $A_0$ is singular or $A_0^{1/\alpha}$ does not exists then (5.32) is replaced by a corresponding less than or equal to relation of norms with the replacements $A_0 \to \|A_0\|$, $A_0^{-1} \to \|A_0\|^{-1}$ and $e^{A_0 t} \to \|e^{A_0 t}\|$.

**(iii)** Assume that $J_{A_0} = J_{A_{0d}} + \tilde{J}_{A_0}$ is the canonical real form of $A_0$ ( in particular, its Jordan form if all the eigenvalues are real) with $J_{A_{0d}}$ being diagonal and $\tilde{J}_{A_0}$ being off-diagonal such that the above decomposition is unique with $A_0 = T^{-1} J_{A_0} T$ where T is a unique non-singular transformation matrix. Then, the Caputo fractional differential system (3.1) is globally Lyapunov stable independently of $A_0^{1/\alpha}$ to exist or not by replacing $\mu_2\left(A_0^{1/\alpha}\right) \to \mu_2\left(J_{A_0}\right)$ in (5.30) by

$$\left\|\left(\frac{1}{\beta_0} T^{-1}\tilde{J}_{A_0} T, \frac{1}{\beta_1} T^{-1} A_1 T, \frac{1}{\beta_2} T^{-1} A_2 T, \cdots, \frac{1}{\beta_p} T^{-1} A_p T\right)\right\|_2 \leq \left|\mu_2\left(J_{0d}\right)\right|^{1/\alpha} \tag{5.33}$$

with $\mu_2\left(J_{0d}\right) \leq 0$ for some set of numbers $\beta_i \in \mathbf{R}_+$ $\left(i \in \overline{p} \cup \{0\}\right)$ satisfying $\sum_{i=0}^{p} \beta_i^2 = 1$. The fractional system is globally asymptotically Lyapunov stable for one such a set of real numbers if $\mu_2\left(J_{0d}\right) < 0$, and

$$\left\|\left(\frac{1}{\beta_0} T^{-1}\tilde{J}_{A_0} T, \frac{1}{\beta_1} T^{-1} A_1 T, \frac{1}{\beta_2} T^{-1} A_2 T, \cdots, \frac{1}{\beta_p} T^{-1} A_p T\right)\right\|_2 < \left|\mu_2\left(J_{0d}\right)\right|^{1/\alpha} \tag{5.34}$$

**Proof**: It turns out that $x_\alpha(t)$ is bounded for all time so that (3.1) is globally Lyapunov stable if $\|\Psi_{\alpha j 0}(t)\|$ is bounded; $\forall j \in \overline{k-1} \cup \{0\}$ for all $t \in \mathbf{R}_{0+}$ for any bounded functions of initial conditions $\varphi_j : [-h, 0] \to \mathbf{R}^n$; $\forall j \in \overline{k-1} \cup \{0\}$ with $\varphi_j(0) = x_j(0) = x_{j0}$. If, in addition, $\|\Psi_{\alpha j 0}(t)\| \to 0$ as $t \to \infty$ then $x_\alpha(t) \to 0$ as $t \to \infty$ so that (3.1) is globally asymptotically Lyapunov stable and the solution (5.27) is bounded for all time. Thus:

If $k = \alpha = 1$ (the particular standard bon- fractional case) then (3.1) is globally Lyapunov stable if $\|\Psi_{100}(t)\|$ is bounded for all $t \in \mathbf{R}_{0+}$. A sufficient condition independent of the delays is that (5.25) holds requiring trivially for the $\ell_2$- matrix measure of $A_0$ to fulfil $\mu_2(A_0) := \frac{1}{2}\lambda_{max}(A_0 + A_0^T) \leq 0$, where the for some $\beta_i \in \mathbf{R}_+$ $(i \in \bar{p})$ subject to $\sum_{i=1}^{p}\beta_i^2 = 1$, [6]. Eq. (5.26) follows from (5.16) after inspection for N=1 and it is bounded as $t \to \infty$ e since otherwise the global stability property (5.25) would fail contradicting its sufficient condition for $j+1 = k = \alpha = 1$. Eq.(5.27) follows from (5.20) for $j+1 = k = \alpha = N = 1$ irrespective of $A_0$ being singular or non-singular and $A_0^{1/\alpha}$ to exist or not. Eq. (5.28) follows from (5.26) since $\mu_2(A_0) < 0$ implies that $A_0$ is a stability matrix then $Re(\lambda) < 0$; $\forall \lambda \in \sigma(A_0)$ and, furthermore, $\Psi_{100}(t) \to 0$, and the unforced solution $x_\alpha(t) \to 0$, as $t \to \infty$ from the strict inequality guaranteeing global asymptotic stability independent of the delays, namely, $\left\|\frac{1}{\beta_1}A_1, \frac{1}{\beta_2}A_2, \cdots, \frac{1}{\beta_p}A_p\right\|_2 < |\mu_2(A_0)|$. Property (i) has been proven. Property (ii) has a similar proof for $\alpha \in (0,1]$, k=1 by replacing $A_0 \to A_0^{1/\alpha}$. Property (iii) follows by using the matrix similarity transformation $A_0 = T^{-1}J_{A_0}T = T^{-1}(J_{A_{0d}} + \tilde{J}_{A_0})T$ and using the homogeneous transformed Caputo fractional differential system from (3.1):

$$\left({}^{C}D_{0+}^{\alpha}z\right)(t) = \left({}^{C}D_{0+}^{\alpha}Tx\right)(t) = \sum_{i=0}^{p}A_i Tx(t-h_i) \Leftrightarrow$$

$$\left({}^{C}D_{0+}^{\alpha}x\right)(t) = \sum_{i=0}^{p}T^{-1}A_i Tx(t-h_i) = T^{-1}A_0 Tx(t) + \sum_{i=1}^{p}T^{-1}A_i Tx(t-h_i)$$

$$= T^{-1}J_{A_{0d}}Tx(t) + \sum_{i=0}^{p}T^{-1}\bar{A}_i Tx(t-h_i) \qquad (5.35)$$

where $z(t) = Tx(t)$; $\forall t \in \mathbf{R}_{0+}$, $h_0 = 0$ plays the role of an additional delay. $\bar{A}_0 = \tilde{J}_{A_0}$ and $\bar{A}_i = A_i$ $(i \in \bar{p})$ by noting also that since $\left(J_{A_{0d}} + J_{A_{0d}}^*\right)$ is diagonal with real eigenvalues by construction, one has:

$$\left|\mu_2\left(J_{A_{0d}}^{1/\alpha}\right)\right| = \left|\frac{1}{2}\lambda_{max}\left(J_{A_{0d}}^{1/\alpha} + J_{A_{0d}}^{1/\alpha *}\right)\right| = \left|\lambda_{max}\left(J_{A_{0d}}^{1/\alpha}\right)\right|$$

$$= \left|Re\,\lambda_{max}\left(J_{A_{0d}}^{1/\alpha}\right)\right| = \left|Re\,\lambda_{max}^{1/\alpha}\left(J_{A_{0d}}\right)\right| = \left|Re\,\lambda_{max}^{1/\alpha}(A_{0d})\right| = \left|\mu_2(J_{0d})\right|^{1/\alpha} \quad (5.36)$$

Then, the proof is similar to that of the related part of Property (ii). $\square$

*Remark 5.2.* Note that a similar expressions to (5.32) applies to guarantee global asymptotic stability for $\alpha \in (0,1]$ in Theorem 5.1(iii) by replacing $A_0 \to T^{-1}J_{A_{0d}}T$ and $A_i \to T^{-1}\bar{A}_i T$ with $\overline{A}_i$ $(i \in \overline{p} \cup \{0\})$ defined in the proof of Theorem 5.1(iii). Theorem 5.1 establishes that for any stability matrix $A_0$, the asymptotic stability condition of sufficient type is as follows:

$$\left\| \left( \frac{1}{\beta_0} T^{-1} \widetilde{J}_{A_0} T, \frac{1}{\beta_1} T^{-1} A_1 T, \frac{1}{\beta_2} T^{-1} A_2 T, \cdots, \frac{1}{\beta_p} T^{-1} A_p T \right) \right\|_2 < \left| \mu_2 (J_{0d}) \right|^{1/\alpha} \qquad (5.37)$$

provided that $\mu_2(J_{0d}) = Re\, \lambda_{max}(A_0) < 0$ extends from $\alpha = \alpha_0 \leq 1$, (in particular, from the standard non-fractional differential system $\alpha = \alpha_0 = 1$) to any $\alpha \in (0, \alpha_0]$

$$\left\| \left( \frac{1}{\beta_0} T^{-1} \widetilde{J}_{A_0} T, \frac{1}{\beta_1} T^{-1} A_1 T, \frac{1}{\beta_2} T^{-1} A_2 T, \cdots, \frac{1}{\beta_p} T^{-1} A_p T \right) \right\|_2 < \left| \mu_2 (J_{0d}) \right|^{1/\alpha_0} \leq \left| \mu_2 (J_{0d}) \right|^{1/\alpha}$$

$$; \forall \alpha \in (0, \alpha_0] \qquad (5.38)$$

Note that the global Lyapunov´s stability conditions (5.30) and (5.33) with nonpositive measures $\mu_2(J_{0d})$ being eventually zero of the corresponding matrices of the unforced fractional dynamic system does not imply the boundedness of the solutions of the system for any admissible forcing bounded control. However, under strict inequalities (5.31) or (5.34) and negative related matrix measures $\mu_2(J_{0d})$, i.e. if asymptotic stability holds, the forced solutions for any bounded controls are guaranteed to be uniformly bounded. □

It follows after inspecting the solution (3.8), subject to (3.9)-(3.10), and the expressions (5.22)-(5.23) that the stability properties for arbitrary admissible initial conditions or admissible bounded controls are lost in general if $\alpha \geq 2$. However, it turns out that the boundedness of the solutions can be obtained by zeroing some of the functions of initial conditions. Note, in particular, from (5.22)-(5.23) that $\varphi_j$ is required to be identically zero on its definition domain for $\overline{k-1} \cup \{0\} \ni j < \alpha - 1$ $(\alpha \geq 2)$ in order that the $\Gamma$-functions be positive (note that $\Gamma(x)$ is discontinuous at zero with an asymptote to $-\infty$ as $x \to 0^-$). This observation combined with Theorem 5.1 leads to the following direct result which is not a global stability result:

**Theorem 5.2**. Assume that $\alpha \geq 2$ and the constraint (5.32) holds with negative matrix measure $\mu_2(J_{0d})$. Assume also that $\varphi_j : [-h, 0] \to \mathbf{R}^n$ are any admissible functions of initial conditions for $\overline{k-1} \cup \{0\} \ni j \geq \alpha - 1$ while they are identically zero if $\overline{k-1} \cup \{0\} \ni j < \alpha - 1$. Then, the unforced solutions are uniformly bounded for all time independent of the delays. Also, the total solutions for admissible bounded controls are also bounded for all time independent of the delays. □

The stability of positive or nonnegative solutions is of direct characterization by combining the positivity conditions of the above section with the stability analysis of this section.

**ACKNOWLEDGMENTS**

The author is grateful to the Spanish Ministry of Education by its partial support of this work through Grant DPI2009-07197. He is also grateful to the Basque Government by its support through Grants GIC07143-IT-269-07and SAIOTEK S-PE08UN15.

**REFERENCES**

[1] A. K. Kilbas, H.M. Srivastava and J. J. Trujillo, *Theory and Applications of Fractional Differential Equations*, North- Holland Mathematics Studies, 204, Jan van Mill Editor, Elsevier, Amsterdam, 2006.

[2] Z. Odibat, "Approximations of fractional integrals and Caputo fractional derivatives", *Applied Mathematics and Computation*, Vol. 178, pp. 527-533, 2006.

[3] Y. Luchko and R. Gorenflo, "An operational method for solving fractional differential equations with the Caputo derivatives", *Acta Mathematica Vietnamica*, Vol. 24, No. 2, pp. 207-233, 1999.

[4] Y. Luchko and H.M. Srivastava, "The exact solution of certain differential equations of fractional order by using operational calculus", *Computers and Mathematics with Applications*, Vol. 29, No. 8. pp. 73-85, 1995.

[5] R.C. Soni and D. Singh, "Certain fractional derivative formulae involving the product of a general class of polynomials and the multivariable H- function" , *Proc. Indian Acad. Sci.( Math. Sci.)*, Vol. 112, No. 4, pp. 551- 562, 2002.

[6] R.K. Raina, "A note on the fractional derivatives of a general system of polynomials" , *Indian Journal of Pure and Applied Mathematics*, Vol. 16, No. 7 , pp 770-774, 1985.

[7] A. Babakhani and V.D.Gejji, "On calculus of local fractional derivatives", *Journal of Applied Mathematics and Applications*, Vol. 270, pp. 66-79, 2002.

[8] A. Ashyralyev, "A note on fractional derivatives and fractional powers of operators", *Journal of Applied Mathematics and Applications*, Vol. 357, pp. 2342-236, 2009.

[9]B. Baeumer, M. M. Meerschaert and J. Mortensen, "Space-time fractional derivative operators", *Proceedings of the American Mathematical Society*, Vol. 133, No. 8 , pp. 2273-2282, 2005.

[10] B. Ahmad and S. Sivasundaram, " On four-point nonlocal boundary value problems of nonlinear integro-differential equations of fractional order", *Applied Mathematics and Computation* , Vol. 217, No. 2, pp. 480-487, 2010.

[11] F. Riewe, "Merchanics with fractional derivatives", *Physical Review E*, Vol. 55, No. 3 pp. 3581-3592, 1997.

[12] A.M.A. El- Sayed and M. Gaber, "On the finite Caputo and finite Riesz derivatives", *Electronic Journal of Theoretic Physics*, EJTP 3, no. 12, pp. 81-95, 2006.

[13] R. Almeida, A. B. Malinowska, D.F. M. Torres, "A fractional calculus of variations for multiple integrals with applications to vibrating string", *arXiv:1001.2722v1[math.OC] 15 Jan 2010*. Also in *Journal of Mathematical Physics* (in press).

[14] I. Schaffer and S. Kempfle, "Impulse responses of fractional damped systems", *Nonlinear Dynamics*, Vol. 38, pp. 61-68, 2004.

[15] F.J. Molz III, G.J. Fix III and S. Lu, "A physical interpretation of the fractional derivative in Levy diffusion", *Applied Mathematics Letters*, Vol. 15, pp. 907-911, 2002.

[16] M. IIic, IW Turner, F. Liu and V. Anh, "Analytical and numerical simulation of a one-dimensional fractional-in-space diffusion equation in a composite medium", *Applied Mathematics and Computation*, Vol. 216, No. 8 , pp.2248-2262, 2010.

[17] M. D. Ortigueira, "On the initial conditions in continuous- time fractional linear systems", *Signal Processing* ,Vol. 83 , No. 3 , pp 2301-2309, 2003.

[18] J. Lancis, T. Szoplik, E. Tajahuerce, V. Climent and M. Fernandez- Alonso, "Fractional derivative Fourier plane filter for phase – change visualization", *Applied Optics*, Vol. 36 , No. 29, pp. 7461-7464, 1997.

[19] E. Scalas, R. Gorenflo and F. Mainardi, "Fractional calculus and continuous-time finance", *arXiv:cond-mat/ 0001120v1 [cond-mat.dis-nn] 15 Jan 2000*.

[20] S. Das, *Functional Fractional Calculus for System Identification and Controls*, Springer- Verlag, Berlin 2008.

[21] J.A. Tenreiro Machado, Manuel F. Silva, Ramiro S. Barbosa, Isabel S. Jesus, Cecilia M. Reis, Maria G. Marcos and Alexandra F. Galhano, " Some applications of fractional calculus in Engineering", *Mathematical Problems in Engineering*, Vol. 2010, Article ID 639801, 34 pages, doi: 10.1155/2010/639801.

[22] J. Sabatier, O.P. Agrawal and J.A. Tenreiro Machado Eds., *Advances in Fractional Calculus: Theoretical Developments and Applications in Physics and Engineering*, Springer- Verlag, 2007.

[23] L. A. Zadeh and C.A. Desoer, *Linear Systems Theory: The State Space Approach*, McGraw- Hill, New York, 1963.

[24] M. DelaSen, " Application of the non-periodic sampling to the identifiability and model-matching in dynamic systems", *International Journal of Systems Science*, Vol. 14, No. 4, pp. 367-383, 1983.

[25] M. De la Sen and NS Luo, "On the uniform stability of a wide class of time-delay systems", *Journal of Mathematical Analysis and Applications,* vol. 289, No. 2, pp. 456-476, 2004.

[26] M. De la Sen, "Stability of impulsive time-varying systems and compactness of the operators mapping the input space into the state and output space", *Journal of Mathematical Analysis and Applications,* vol. 321, No. 2, pp. 621-650, 2006.

[27] J. Chen, D. Xu and B. Shafai, "On sufficient conditions for stability independent of delay", *IEEE Transactions on Automatic Control*, Vol. 40, pp. 1675-1680, 1995.

[28] M. De la Sen , "On the reachability and controllability of positive linear time-invariant dynamic systems with internal and external incommensurate point delays", *Rocky Mountain Journal of Mathematics*, Vol. 40, No. 1, pp. 177-207, 2020.

[29] M. De la Sen, "A method for general design positive real functions", *IEEE Transactions on Circuits and Systems I – Fundamental Theory and Applications*, Vol. 45, No. 7, pp. 764-769, 1998.

[30] M. De la Sen, "Preserving positive realness through discretization" , *Positivity*, Vol. 6, No. 1, pp. 31-45, 2002.

[31] M. De la Sen, " On positivity of singular regular linear time-delay time-invariant systems subject to multiple internal and external incommensurate point delays, *Applied Mathematics and Computation*, Vol. 190, No. 1, pp. 382-401, 2007.